\documentclass[10pt]{article}
\usepackage{amsfonts}
\usepackage{hyperref}
\usepackage{mathrsfs,enumerate,amssymb,amsmath,color,lineno}
\usepackage{indentfirst}
\usepackage{amsmath}
\usepackage{amssymb}
\usepackage[amsmath,thmmarks]{ntheorem}
\usepackage{cite}
\usepackage{graphicx}
\usepackage{tikz}

\newtheorem{definition}{\bf Definition}[section]
\newtheorem{lemma}{\bf Lemma}[section]
\newtheorem{theorem}{\bf Theorem}[section]
\newtheorem{remark}{\bf Remark}[section]
\newtheorem{corollary}{\bf Corollary}[section]
\newtheorem{example}{\bf Example}[section]

\newtheorem{proposition}{\bf Proposition}[section]

\begin{document}
\setcounter{page}{1}

\title{{\textbf{Associativity of a class of two-place functions and its consequences for classes of triangular norms
%Characterizations of associativity of two-place functions generated by strictly increasing functions
}}\thanks {Supported by
the National Natural Science Foundation of China (Nos.12071325,12471440)}}
\author{Yun-Mao Zhang\footnote{\emph{E-mail address}: 2115634219@qq.com}, Xue-ping Wang\footnote{Corresponding author. xpwang1@hotmail.com; fax: +86-28-84761502 (id:0000-0003-4376-9859)},\\
\emph{School of Mathematical Sciences, Sichuan Normal University,}\\
\emph{Chengdu 610066, Sichuan, People's Republic of China}}

\newcommand{\pp}[2]{\frac{\partial #1}{\partial #2}}
\date{}
\maketitle
\begin{quote}
{\bf Abstract} This article characterizes the associativity of two-place functions $T: [0,1]^2\rightarrow [0,1]$ defined by $T(x,y)=f^{(-1)}(F(f(x),f(y)))$ where $F:[0,1]^2\rightarrow[0,1]$ is a triangular norm (even a triangular subnorm), $f: [0,1]\rightarrow [0,1]$ is a strictly increasing function and $f^{(-1)}:[0,1]\rightarrow[0,1]$ is the pseudo-inverse of $f$. We prove that the associativity of functions $T$ only depends on the range of $f$, which is used to give a sufficient and necessary condition for the function $T$ being associative when the triangular norm $F$ is an ordinal sum of triangular norms and an ordinal sum of triangular subnorms in the sense of A. H. Clifford, respectively. These results finally are applied for describing classes of triangular norms generated by strictly increasing functions.

%{\textbf{\emph{Mathematics Subject Classification}}}: 20M05
\end{quote}

\section{Introduction}
%Triangular norms (t-norms for short) have been introduced in the framework of statistical  metric spaces for proper generalization of the triangle inequality from metric spaces to statistical metric spaces \cite{CA2006, EP2000, BS1960}. In the context of fuzzy logic where it is a many-valued propositional logic with a continuum of truth values modelled by unit intervals, t-norms and triangular conorms (t-conorms for short) model the (semantic) interpretations of the logical conjunction and disjunction, respectively. The condition of left-continuity of t-norms is a frequently cited property and plays a central role in all the fields that use t-norms.  For example, comparing with all the other approaches, the crucial role of logics with implications defined as the residuum of the conjunction is stressed, and it is well known that the left-continuity of a t-norm and the definability of the (t-norm based) residual implication are equivalent. From the usually cited paper of Ling \cite{EP1999}, continuous t-norms have become well understood and have been used in several applications.

It is well known that the constructing methods of triangular norms (t-norms for short) play a central role in the theory of t-norms. There are many methods in the current literature concerning the construction of t-norms, see e.g., \cite{CA2006,AC1954,SJ1999,SJ2002,Klement1999,EP2000,EP2002,CL1965,AM2004,YO2007,YO2008,BS1983,PV1998,PV2005,PV2008,PV2010,PV2013,PV2020,PV2022,DZ2005}. For example, based on known t-norms $T: [0,1]^2\rightarrow [0,1]$, monotone functions $f: [0,1]\rightarrow [0,1]$ and their pseudo-inverses $f^{(-1)}$, Klement, Mesiar and Pap \cite{Klement1999,EP2000} gave a sufficient condition for the following function $T_{[f]}: [0,1]^2\rightarrow [0,1]$ given by
 \begin{equation*}
T_{[f]}(x,y)=\left\{
  \begin{array}{ll}
    \min\{x,y\} & \mbox{if }\max\{x,y\}=1, \\
  f^{(-1)}(T(f(x),f(y))) & \mbox{otherwise,}
  \end{array}
\right.
\end{equation*}
being a t-norm. A lot of specific examples for this construction are provided and their properties are studied. Moreover, Klement, Mesiar and Pap \cite{EP2000} used the known operations, the usual addition and multiplication of real numbers, to construct t-norms in terms of one-place functions. They defined an additive generator $f :[0,1]\rightarrow [0,\infty]$ of a t-norm $T$ as a strictly decreasing function which is also right-continuous at $0$ with $f (1) = 0$, such that for all $(x,y)\in [0,1]^2$
\begin{equation}
\label{eq:1.7}
f(x)+f(y)\in \mbox{Ran}(f)\cup [f(0),\infty],
\end{equation}

\begin{equation}
\label{eq:1.8}
 T(x,y)= f^{(-1)}(f(x)+f(y))
\end{equation} and thought that one can generalize the concept of additive generator of t-norms in the sense that a strictly decreasing function $f :[0,1]\rightarrow [0,\infty]$ with $f (1) = 0$ can generate a t-norm $T$ via Eq.\eqref{eq:1.8}, i.e., such that for all $(x,y)\in [0,1]^2$ we have $ T(x,y)= f^{(-1)}(f(x)+f(y))$, without satisfying Eq.\eqref{eq:1.7}. After that, scholars all over the world provided a large number of results concerning generalized additive generators of associative functions \cite{Klement1999,PV1998,PV2005,PV2008,PV2010,PV2013,PV2020,PV2022}. For example, Vicen\'{\i}k \cite{PV2005} showed some characterizations for the function $T$ given by Eq.\eqref{eq:1.8} being associative. Vicen\'{\i}k \cite{PV2008} further presented some necessary and sufficient conditions for the function $T$ given by Eq.\eqref{eq:1.8} to be a border-continuous t-norm. Vicen\'{\i}k \cite{PV2013} also explored the generated border-continuous t-norms and verified what properties generated border-continuous t-norms can (or cannot) posses. In particular, Zhang and Wang \cite{YM2024} recently investigated the associativity of function $T:[0,1]^2\rightarrow [0,1]$ defined by \begin{equation}\label{eq1}
T(x,y)= f^{(-1)}(F(f(x),f(y)))
\end{equation}
where $f^{(-1)}$ is the pseudo-inverse of monotone right continuous function $f: [0,1]\rightarrow [0,1]$, and $F: [0,1]^2\rightarrow [0,1]$ is a monotone associative function with neutral element in $[0,1]$. Another interesting and powerful tool of the constructing methods for t-norms is an ordinal sum of t-norms which was introduced by \cite{BS1963,BS1983}. It is shown that the ordinal sum theorem for a t-norm over unit interval remains valid if we replace the requirement that the summands are a family of t-norms by they are a family of t-subnorms (with some additional mild properties if necessary) \cite{SJ2002}. Furthermore, Klement, Mesiar and Pap \cite{EP2002} had shown that the generalized ordinal sum theorem for a t-norm %the ordinal sum theorem for a t-norm in \cite{SJ2002}
is the most general way to obtain a t-norm as an ordinal sum of semigroups and all t-norms over unit interval can be naturally divided into three classes, i.e., $\{T_M\}$, the class of ordinally irreducible t-norms and the class of t-norms which are different from $\{T_M\}$ and are not ordinally irreducible. %Therefore, an important unsolved problem is: what are the characterizations of such classes of t-norms in the sense of Klement, Mesiar and Pap \cite{EP2002}? This article will positively answer the problem through the two-place function $T$ given by Eq.~\eqref{eq1}.

However, the characterizations of the associativity of generated functions given by Eq.~\eqref{eq1} are still thoroughly unsolved problems. Motivated by the work of Vicen\'{\i}k \cite{PV2005} together with Zhang and Wang \cite{YM2024}, this article pays attention to the associativity of two-place functions given by Eq.~\eqref{eq1} when $f$ is a strictly increasing function and $F$ is a t-norm (even a t-subnorm).
%Since $F$ is a t-norm that can be represented as an ordinal sum of t-subnorms by \cite{SJ2002,EP2002}, a sufficient and necessary condition for the function $T$ defined by Eq.~\eqref{eq1} being t-norm is further given, and its consequence is applied for describing the classes of triangular norms. The work of this paper has potential applications in constructing t-norms, classifying associative functions and simplifying the related calculations.

The rest of this article is organized as follows. In Section 2, we recall some known concepts and results of t-norms and properties of strictly monotone functions, respectively. The work of Section 3 is developed when $f$ is a strictly increasing function and $F$ is a t-norm, in which some necessary and sufficient conditions for the associativity of function $T$ via Eq.\eqref{eq1} are given in terms of the range of the $f$. In Sections 4, we use the results of Section 3 to investigate the associativity of function $T$ given by Eq.~\eqref{eq1} when the t-norm $F$ is an ordinal sum of t-norms in the sense of A. H. Clifford. In Section 5, we first describe the associativity of function $T$ given by Eq.~\eqref{eq1} when $F$ is a t-subnorm, and then explore the associativity of function $T$ given by Eq.~\eqref{eq1} when the t-norm $F$ is an ordinal sum of t-subnorms in the sense of A. H. Clifford. In Section 6, three necessary and sufficient conditions are presented for the function $T$ given by Eq.~\eqref{eq1} being divided into three classes: $\{T_M\}$, ordinally irreducible t-norm and t-norm with a non-trivial ordinal sum representation and $T\neq T_M$, respectively. A conclusion is drawn in Section 7.

\section{Preliminaries}

In this section, we recall some elementary concepts and results.

Let $[p,q]\subseteq [-\infty,\infty]$ with $p\leq q$. Then by convention, $\sup \emptyset=p$ and $\inf \emptyset=q$.
\begin{definition}[\cite{EP2000,PV2005}]\label{def2.1}
\emph{Let $p, q, s, t\in [-\infty, \infty]$ with $p<q, s<t$ and $f:[p,q]\rightarrow[s,t]$ be a non-decreasing (resp. non-increasing) function. Then the function $f^{(-1)}:[s,t]\rightarrow[p,q]$ defined by
\begin{equation*}
f^{(-1)}(y)=\sup\{x\in [p,q]\mid f(x)<y\}\,(\mbox{resp. }f^{(-1)}(y)=\sup\{x\in [p,q]\mid f(x)>y\})
\end{equation*}
is called the \emph{pseudo-inverse} of a non-decreasing (resp. non-increasing) function $f$. In particular, $f^{-1}(y)=\{x\in [p,q]\mid f(x)=y\}$.}
\end{definition}

The pseudo-inverse of non-decreasing (resp. non-increasing) function is also non-decreasing (resp. non-increasing). Moreover, if $f$ is strictly monotone, then its pseudo-inverse $f^{(-1)}$ is strictly monotone on $\mbox{Ran}(f)$ where $\mbox{Ran}(f)=\{f(x)\mid x\in[p,q]\}$, and continuous on $[0,1]$. In particular, if $f$ is strictly monotone, then $f^{(-1)}(f(x))=x$ for all $x\in[p,q]$, and for all $x\in[s,t]$, $f(f^{(-1)}(x))=x$ if and only if $x\in \mbox{Ran}(f)$.

Let $M\subseteq [0,1]$. A point $x\in [0,1]$ is said to be an \emph{accumulation point of $M$ from the left} if there exists a strictly increasing sequence $\{x_{n}\}_{n\in N}$ of points $x_{n}\in M$ such that $\lim_{n\rightarrow \infty}x_{n}=x$. A point $x\in [0,1]$ is said to be an \emph{accumulation point of $M$ from the right} if there exists a strictly decreasing sequence $\{x_{n}\}_{n\in N}$ of points $x_{n}\in M$ such that $\lim_{n\rightarrow \infty}x_{n}=x$. A point $x\in [0,1]$ is said to be an \emph{accumulation point} of $M$ if $x$  is an accumulation point of $M$ from the left or from the right. Denote the set of all accumulation points of $M$ from the left (resp. right) by $ACC_{-}(M)$ (resp. $ACC_{+}(M)$) and the set of all accumulation points of $M$ from the left and right by $ACC_{0}(M)$, respectively. It is easy to see that $ACC_{-}(M)\cap ACC_{+}(M)=ACC_{0}(M)$.

For a monotone function $f:[0,1]\rightarrow [0,1]$, let $f(a^-)=\lim_ {x \rightarrow a^{- }} f(x)$  for each $a\in (0, 1]$ and $f(a^+)=\lim_ {x \rightarrow a^{+}} f(x)$ for each $a\in [0, 1)$. In particular, let $f(0^-)=f(0)$, and $f(1^+)=f(1)$. A strictly increasing function $f:[0,1]\rightarrow[0,1]$ is continuous at $x\in[0, 1]$ if and only if $f(x^-)=f(x)=f(x^+)$. Obviously, a strictly increasing function $f$ is continuous at $x\in(0, 1)$ if and only if $f(x)\in ACC_{0}(\mbox{Ran}(f))$. The set of all points at which a strictly increasing function $f$ is not continuous is countable. If a strictly increasing function $f$ is not continuous at $x\in(0, 1)$, then $f(x^-)\in ACC_{-}(\mbox{Ran}(f))$ and $f(x^+)\in ACC_{+}(\mbox{Ran}(f))$. For a strictly increasing function $f$, between arbitrary two points $x, y\in[0,1]$ with $x < y$ there are infinitely many points at which $f$ is continuous, or equivalently, between arbitrary two points $x, y\in \mbox{Ran}(f)$ with $x < y$, $ACC_{0}(\mbox{Ran}(f))$ has an infinite elements. Therefore, for a strictly increasing function $f$, $\mbox{Ran}(f)\cap[x, y]$ with $0\leq x\leq y\leq 1$ is an empty set, or a singleton set or an infinite set.

In complete analogy to \cite{PV2013}, we have the following properties of the pseudo-inverse $f^{(-1)}$ of a strictly increasing function $f:[0,1]\rightarrow[0,1]$  which are often used latter.
\renewcommand{\labelenumi}{(\roman{enumi})}
\begin{enumerate}
\item If $x>f(t^+)$ and $t\in[0, 1)$, then $t< f^{(-1)}(x)$. If $x\leq f(t^+)$ and $t\in[0, 1)$, then $f^{(-1)}(x)\leq t$. If $x<f(t^-)$ and $t\in(0, 1]$, then $f^{(-1)}(x)<t$. If $x\geq f(t^-)$ and $t\in(0, 1]$, then $f^{(-1)}(x)\geq t$.
\item If $f(t^-)\leq x\leq f(t^+)$ with $t\in(0, 1)$, then $t=f^{(-1)}(x)$. If $x\leq f(0^+)$ then $f^{(-1)}(x)=0$. If $x\geq f(1^-)$ then $f^{(-1)}(x)=1$.
\item Let $0\leq x\leq y\leq 1$. Then $f^{(-1)}(x)=f^{(-1)}(y)$ if and only if $\mbox{Ran}(f)\cap[x, y]$ is at most a one-element set.
\end{enumerate}

\begin{definition}[\cite{EP2000}]
\emph{A \emph{t-norm} is a binary operator $T:[0, 1]^2\rightarrow [0, 1]$ such that for all $x, y, z\in[0, 1]$ the following conditions are satisfied:}

$(T1)$  $T(x,y)=T(y,x)$,

$(T2)$  $T(T(x,y),z)=T(x,T(y,z))$,

\emph{$(T3)$  $T(x,y)\leq T(x,z)$ whenever $y\leq z$,}

$(T4)$  $T(x,1)=T(1,x)=x$.
\end{definition}

A binary operator $T:[0, 1]^2\rightarrow [0, 1]$ is called \emph{t-subnorm} if it satisfies $(T1), (T2), (T3)$, and $T(x,y)\leq \min\{x,y\}$ for all $x,y\in [0, 1]$. Further, it is called a \emph{proper t-subnorm} if it is not a t-norm.
\begin{definition}[\cite{EP2000}]\label{def4.1}\emph{(Ordinal Sum Theorem for t-norms). Let $A\neq\emptyset$ be a totally ordered set, $\{F_\alpha \mid \alpha\in A\}$ be a family of t-norms and $\{(a_\alpha,b_\alpha)\}_{\alpha\in A}$ be a family of pairwise disjoint open subintervals of $[0, 1]$. Then \begin{equation*}
F(x,y)=\left\{
  \begin{array}{ll}
    a_\alpha+(b_\alpha-a_\alpha)F_\alpha(\frac{x-a_\alpha}{b_\alpha-a_\alpha},\frac{y-a_\alpha}{b_\alpha-a_\alpha}) & \hbox{if }(x,y)\in[a_\alpha,b_\alpha]^{2}, \\
    \min\{x,y\} & \hbox{otherwise}
  \end{array}
\right.
\end{equation*} is also a t-norm. The t-norm $F$ is called an \emph{ordinal sum} of $\{F_\alpha \mid \alpha\in A\}$, and is written as $F=(\langle a_\alpha,b_\alpha, F_\alpha\rangle)_{\alpha\in A}$. Further, $F_\alpha$ and $(a_\alpha,b_\alpha)$ are referred to as a \emph{summand} and an \emph{open generating interval} of the ordinal sum, respectively.}
\end{definition}
\section{Characterizations of the associativity of function $T$ given by \eqref{eq1}}
In this section, we first define an operation $\otimes$ on the range Ran$(f)$ of a strictly increasing function $f:[0,1]\rightarrow [0,1]$ and then characterizes what properties of Ran$(f)$ are equivalent to the associativity of operation $\otimes$. For the completeness and readability of this
article, we have provided some necessary proofs even if all the concepts and results of this section are originated from \cite{PV2010,YM2024}.

In what follows, we always suppose that $f:[0,1]\rightarrow [0,1]$ is a strictly increasing function.

Denoted by $A\setminus B=\{x\in A\mid x\notin B\}$ for two sets $A$ and $B$. Let $M\subseteq[0,1]$ and
$$\mathcal{A}=\{M \mid \mbox{there is a strictly increasing function } f:[0,1]\rightarrow[0,1]\mbox{ such that }\mbox{Ran}(f)=M\}.$$
Note that there exist many strictly increasing functions with the same range.

It is well known that a strictly decreasing function from $[0,1]\rightarrow[0,\infty]$ can be transformed into a strictly increasing function from $[0,1]\rightarrow[0,1]$ by a simple exponential transformation. Therefore, in complete analogy to \cite{PV2010} and \cite{YM2024}, we can give the following two definitions.
\begin{definition}\label{3.0}
\emph{Let $M\in \mathcal{A}$. A pair $(\mathcal{S},\mathcal{C})$ is said to be \emph{associated} with $M\neq [0,1]$ if $\mathcal{S}=\{[b_k, d_k] \subseteq [0,1]\mid k\in K\}$ is a non-empty countable system of pairwise disjoint intervals of a positive length and $\mathcal{C}=\{c_k\in [0,1]\mid k\in K\}$ is a non-empty countable set such that $[b_k, d_k] \cap \mathcal{C}=\{c_k\}$ for all $k\in K$ and
\begin{equation*}
M= \{c_k\in [0,1]\mid k\in K\}\cup \left([0,1]\setminus \left(\bigcup_{k\in K}[b_k, d_k] \right)\right).
\end{equation*}
A pair $(\mathcal{S},\mathcal{C})$ is said to be \emph{associated} with $M=[0,1]$ if $\mathcal{S}=\{[1,1]\}$ and $\mathcal{C}=\{1\}$. We shortly write $(\mathcal{S},\mathcal{C})=(\{[b_k, d_k]\mid k\in K\}, \{c_k \mid k\in K\})$ instead of $(\mathcal{S},\mathcal{C})=(\{[b_k, d_k]\subseteq [0,1]\mid k\in K\}, \{c_k\in [0,1] \mid k\in K\})$.}
\end{definition}

 Observe that $[0,1]\setminus (M\setminus\mathcal{C})=\bigcup_{k\in K}[b_k, d_k]$ and $[b_k, d_k]\cap M=[b_k, d_k]\cap \mathcal{C}=\{c_k\}$ for all $k\in K$. If $x\in(M\setminus\mathcal{C})\cap(0,1)$ then $x$ is an accumulation point of $M$ from the left and right. Moreover, $x\in[0,1]$ is an accumulation point of $M$ from the left (resp. right) if and only if $x$ is an accumulation point of $M\setminus\mathcal{C}$ from the left (resp. right).

\begin{definition}\label{def3.1}
\emph{Let $M\in \mathcal{A}$. Define a function $G_{M}:[0,1]\rightarrow M $ by
\begin{equation*}
 \{G_{M}(x)|x\in[0,1]\}=M\cap[\sup([0,x]\cap M),\inf([x,1]\cap M)]
\end{equation*}}
\end{definition}

From Definition \ref{def3.1}, it is obvious that $G_{M}(x)=x$ for any $x\in M$ and for any $x\in [0,1]\setminus M$, there exists a $k\in K$ such that $G_{M}(x)=c_k$. Furthermore, we have the following proposition.
\begin{proposition}\label{prop3.1}
Let $M\in \mathcal{A}$ and $(\mathcal{S},\mathcal{C})=(\{[b_k, d_k]\mid k\in K\}, \{c_k \mid k\in K\})$ be associated with $M$. Then for all $x,y\in[0,1]$ and $k\in K$,
\renewcommand{\labelenumi}{(\roman{enumi})}
\begin{enumerate}
\item $G_{M}(x)=x$ if and only if $x\in M$.
\item $G_{M}(x)=c_k$ if and only if $x\in [b_k, d_k]$ and $\{c_k\}=M\cap[b_k, d_k]$.
\item $G_{M}$ is a non-decreasing function.
\end{enumerate}
\end{proposition}

In what follows, we suppose that $F:[0,1]^2\rightarrow [0,1]$ is a t-norm unless otherwise stated.
\begin{definition}\label{def3.2}
\emph{Let $M\in \mathcal{A}$. Define an operation $\otimes:M^2\rightarrow M $ by
\begin{equation*}
x\otimes y=G_{M}(F(x,y)).
\end{equation*}}
\end{definition}

Using Propositions \ref{prop3.1} (i) and (ii) and Definition \ref{def3.2}, we have the following proposition.
\begin{proposition}\label{prop3.2}
Let $M\in \mathcal{A}$ and $(\mathcal{S},\mathcal{C})=(\{[b_k, d_k]\mid k\in K\}, \{c_k \mid k\in K\})$ be associated with $M$. Then for all $x,y\in M$ and $k\in K$,
\renewcommand{\labelenumi}{(\roman{enumi})}
\begin{enumerate}
\item $x\otimes y=F(x,y)$ if and only if $F(x,y)\in M$.
\item $x\otimes y=c_k$ if and only if $F(x,y)\in [b_k, d_k]$ and $\{c_k\}=M\cap[b_k, d_k]$.
\end{enumerate}
\end{proposition}

Propositions \ref{prop3.1} (i) and (ii) also imply the following lemma.
\begin{lemma}\label{lem3.1}
Let $M\in \mathcal{A}$, $(\mathcal{S},\mathcal{C})=(\{[b_k, d_k]\mid k\in K\}, \{c_k \mid k\in K\})$ be associated with $M$ and $f:[0,1]\rightarrow [0,1] $ be a strictly increasing function with $\emph{Ran}(f)=M$. Then $G_{M}(x)=f(f^{(-1)}(x))$ for all $x\in[0,1]$.
\end{lemma}

The next proposition describes the relation between $T$ given by \eqref{eq1} and $\otimes$.
\begin{proposition}\label{prop3.3}
Let $f:[0,1]\rightarrow [0,1]$ be a strictly increasing function with $\emph{Ran}(f)=M$, and $f^{-1}:M\rightarrow [0,1]$ be the inverse function of $f$. Then
 \begin{equation}\label{eq03}
x\otimes y=f(T(f^{-1}(x),f^{-1}(y)))
\end{equation}
for all $x,y\in M$, and
 \begin{equation}\label{eq04}
T(x,y)=f^{-1}(f(x)\otimes f(y))
\end{equation}
for all $x,y\in [0,1]$. Moreover, $T$ is associative if and only if $\otimes$ is associative.
\end{proposition}
\begin{proof} For all $x,y\in M$, it follows from Definition \ref{def3.2} and Lemma \ref{lem3.1} that
\begin{eqnarray*}
x\otimes y&=&G_{M}(F(x,y))\\
&=&f(f^{(-1)}(F(x,y)))\nonumber\\
&=&f(f^{(-1)}(F(f(f^{-1}(x),f(f^{-1}(y))))\\
&=&f(T(f^{-1}(x),(f^{-1}(y))))
\end{eqnarray*}
since $f$ is a strictly increasing function. At the same time, for every $x,y\in M$ there exist two elements $u,v\in [0,1]$ such that $f(u)=x$ and $f(v)=y$, respectively. Thus
\begin{equation*}
f(u)\otimes f(v)=x\otimes y=f(T(f^{-1}(x),(f^{-1}(y))))=f(T(u,v)).
\end{equation*}
This follows that $T(x,y)=f^{-1}(f(x)\otimes f(y))$ for all $x,y\in [0,1]$.

Supposing $T$ is associative, we prove that $\otimes$ is associative, i.e., $(x\otimes y)\otimes z=x\otimes (y\otimes z)$ for all $x,y,z\in M$. It is clear that $f^{-1}(f(T(f^{-1}(x),(f^{-1}(y)))))=T(f^{-1}(x),(f^{-1}(y)))$  for all $x,y\in M$. Now, let $x,y,z\in M$. Then by the associativity of $T$ and \eqref{eq03}, we have
\begin{eqnarray*}
(x\otimes y)\otimes z&=&f(T(f^{-1}(f(T(f^{-1}(x),(f^{-1}(y))))),f^{-1}(z)))\\
&=&f(T(T(f^{-1}(x),(f^{-1}(y))),f^{-1}(z)))\\
&=&f(T(f^{-1}(x),T(f^{-1}(y),f^{-1}(z))))\\
&=&f(T(f^{-1}(x),f^{-1}(f(T(f^{-1}(y),f^{-1}(z))))))\\
&=&f(T(f^{-1}(x),f^{-1}(y\otimes z)))\\
&=&x\otimes (y\otimes z).
\end{eqnarray*}

Next, we prove that $T$ is associative when $\otimes$ is associative. It is clear that $f(f^{-1}(f(x)\otimes f(y)))=f(x)\otimes f(y)$ for all $x,y\in [0,1]$. Let $x,y,z\in [0,1]$. Then by associativity of $\otimes$ and \eqref{eq04}, we have
\begin{eqnarray*}
T(T(x,y),z)&=&f^{-1}(f(f^{-1}(f(x)\otimes f(y)))\otimes f(z))\\
&=&f^{-1}((f(x)\otimes f(y))\otimes f(z))\\
&=&f^{-1}(f(x)\otimes (f(y)\otimes f(z)))\\
&=&f^{-1}(f(x)\otimes (f(f^{-1}(f(y)\otimes f(z)))))\\
&=&f^{-1}(f(x)\otimes f(T(y,z)))\\
&=&T(x,T(y,z)).
\end{eqnarray*}

Therefore, $T$ is associative if and only if $\otimes$ is associative.
\end{proof}

In the sequel, we show a necessary and sufficient condition for the operation $\otimes$ being associative which answers what properties of $M$ are equivalent to the associativity of operation $\otimes$.

Let $M\subseteq [0,1]$. Define $O(M) = \bigcup_{x,y\in M}(\min\{x,y\}, \max\{x,y\})$ whenever $M\neq \emptyset$ (where $(x, x)=\emptyset$), and $O(M)=\emptyset$ whenever $M=\emptyset$. It is obvious that $O(M)=\emptyset$ if and only if $M=\emptyset$ or $M$ contains exact one element.
Let $\emptyset\neq A,B\subseteq [0,1]$ and $c\in[0,1]$. Denote $F(A,c)=\{F(x,c)\mid x\in A\}$, $F(c,A)=\{F(c,x)\mid x\in A\}$ and $F(A,B)=\{F(x,y)\mid x\in A, y\in B\}$.

The following definition is needed.
\begin{definition}\label{def3.3}
\emph{Let $M\in \mathcal{A}$ and $(\mathcal{S},\mathcal{C})=(\{[b_k, d_k]\mid k\in K\}, \{c_k \mid k\in K\})$ be associated with $M$. For all $y\in M$, and $k,l\in K$, define $$M_{k}^{y}=\{x\in M \mid F(x,y)\in [b_k, d_k]\},$$
$$M^{y}=\{x\in M \mid F(x,y)\in M\setminus\mathcal{C} \},$$ $$I_{k}^{y}=O(\{c_k\}\cup F(M_{k}^{y},y)),$$
\begin{equation*}
J_{k,l}^{y}=\begin{cases}
O(F(M_{k}^{y},c_l)\cup F(c_k,M_{l}^{y})), & \hbox{if }\ M_{k}^{y}\neq\emptyset,  M_{l}^{y}\neq\emptyset,\\
\emptyset, & \hbox{otherwise.}
\end{cases}
\end{equation*}
Put $$\mathfrak{T_{1}}(M)=\bigcup_{y\in M}\bigcup_{k\in K}F(I_{k}^{y},M^{y}),$$ $$\mathfrak{T_{2}}(M)=\bigcup_{y\in M}\bigcup_{k,l\in K}J_{k,l}^{y}.$$
The set $\mathfrak{T}(M)=\mathfrak{T_{1}}(M)\cup\mathfrak{T_{2}}(M)$ is called an $F$-\emph{transformation} of $M$.}
\end{definition}

Then we have two key lemmas as below.
\begin{lemma}\label{lem3.2}
Let $M\in \mathcal{A}$ and $(\mathcal{S},\mathcal{C})=(\{[b_k, d_k]\mid k\in K\}, \{c_k \mid k\in K\})$ be associated with $M$. Let $M_1, M_2\subseteq [0,1]$ be two non-empty sets and $c\in[0,1]$.
Then $F(O(M_1\cup M_2),c)\cap (M\setminus\mathcal{C})\neq\emptyset$ if and only if there exist $x\in M_1$ and $y\in M_2$ such that $(\min\{F(x,c),F(y,c)\},\max\{F(x,c),F(y,c)\})\cap (M\setminus\mathcal{C}) \neq\emptyset.$
\end{lemma}
\begin{proof}
Suppose that $F(O(M_1\cup M_2),c)\cap (M\setminus\mathcal{C})\neq\emptyset$. Then there exist two elements $s,t\in M_1\cup M_2$ with $s<t$ such that $(F(s,c),F(t,c))\cap (M\setminus\mathcal{C})\neq\emptyset$ since the t-norm $F$ is non-decreasing.

The following proof is split into three cases.
\renewcommand{\labelenumi}{(\roman{enumi})}
\begin{enumerate}
\item The case that one of $\{s,t\}$ is contained in $M_1$ and the other is contained in $M_2$ is obvious.

\item The case that $s,t\in M_1$ with $s<t$. Choose $w\in M_2$. Then there are three subcases as below.
\renewcommand{\labelenumi}{(\roman{enumi})}
\begin{enumerate}
\item If $w\leq s$ then put $x=t$ and $y =w$.
\item If $w\geq t$ then put $x=s$ and $y =w$.
\item If $s<w<t$, then $(F(s,c),F(w,c))\cap (M\setminus\mathcal{C})\neq\emptyset$ or $(F(w,c),F(t,c))\cap (M\setminus\mathcal{C})\neq\emptyset$. If $(F(s,c),F(w,c))$ $\cap (M\setminus\mathcal{C})\neq\emptyset$ then put $x=s$ and $y =w$. If $(F(w,c),F(t,c))\cap (M\setminus\mathcal{C})\neq\emptyset$ then put $x=t$ and $y =w$.
\end{enumerate}
In any subcase, one can check that $$(\min\{F(x,c),F(y,c)\},\max\{F(x,c),F(y,c)\})\cap (M\setminus\mathcal{C})\neq\emptyset.$$
\item The case that $s,t\in M_2$ is completely analogous.
\end{enumerate}

The converse implication is obvious.
\end{proof}

\begin{lemma}\label{lem3.3}
Let $M\in \mathcal{A}$ and $(\mathcal{S},\mathcal{C})=(\{[b_k, d_k]\mid k\in K\}, \{c_k \mid k\in K\})$ be associated with $M$. Then for any $x,y\in[0,1]$, $G_M(x)\neq G_M(y)$ if and only if $(\min\{x,y\},\max\{x,y\})\cap (M\setminus\mathcal{C})\neq\emptyset$.
\end{lemma}
\begin{proof}If $G_M(x)\neq G_M(y)$, then from Proposition \ref{prop3.1}, $G_M(\min\{x,y\})=s<t=G_M(\max\{x,y\})$ with $s,t\in M$. Obviously $(s,t)\cap (M\setminus\mathcal{C})\neq\emptyset$. Choose $w\in (s,t)\cap (M\setminus\mathcal{C})$. Then from Proposition \ref{prop3.1} and $w=G_M(w)$, we have $G_M(\min\{x,y\})< G_M(w)<G_M(\max\{x,y\})$, thus $\min\{x,y\}< w<\max\{x,y\}$. Hence $(\min\{x,y\},\max\{x,y\})\cap (M\setminus\mathcal{C})\neq\emptyset$.

Conversely, choose $w\in(\min\{x,y\},\max\{x,y\})\cap (M\setminus\mathcal{C})$. Then the point $w$ is an accumulation point of $M$ from the left and right. This follows that $G_M(\min\{x,y\})< G_M(w)<G_M(\max\{x,y\})$, i.e., $G_M(x)\neq G_M(y)$ since the non-decreasing function $G_M$ is strictly increasing on $M$.
\end{proof}

\begin{theorem}\label{them3.1}
Let $M\in \mathcal{A}$, $(\mathcal{S},\mathcal{C})=(\{[b_k, d_k]\mid k\in K\}, \{c_k \mid k\in K\})$ be associated with $M$ and $\mathfrak{T}(M)$ be the $F$-transformation of $M$. Then $(M,\otimes)$ is a semigroup if and only if $\mathfrak{T}(M)\cap (M\setminus\mathcal{C})=\emptyset$.
\end{theorem}
\begin{proof}Let $(\mathcal{S},\mathcal{C})=(\{[b_k, d_k]\mid k\in K\}, \{c_k \mid k\in K\})$ be associated with $M\in \mathcal{A}$. In order to complete the proof, it is enough to prove that the operation $\otimes$ on $M$ is not associative if and only if $\mathfrak{T}(M)\cap (M\setminus\mathcal{C})\neq\emptyset$.

 Suppose that the operation $\otimes$ is not associative, i.e., there exist three elements $x,y,z\in M$ such that $(x\otimes y)\otimes z\neq x\otimes(y\otimes z)$. Then we claim that $F(x,y)\notin M\setminus\mathcal{C}$ or $F(y,z)\notin M\setminus\mathcal{C}$. Otherwise, from the associativity of $F$ and Definition \ref{def3.2}, $F(x,y)\in M\setminus\mathcal{C}$ and $F(y,z)\in M\setminus\mathcal{C}$ would imply $(x\otimes y)\otimes z=G_M(F(F(x,y),z))= G_M(F(x,F(y,z)))= x\otimes(y\otimes z)$, a contradiction. The following proof is split into three cases.

(i) Let $F(x,y)\notin M\setminus\mathcal{C}$ and $F(y,z)\in M\setminus\mathcal{C}$. Then $y\otimes z=F(y,z)$ and there exists a $k\in K$ such that $F(x,y)\in [b_k,d_k]$. Thus $x\otimes y=c_k$ where $\{c_k\}=M\cap[b_k,d_k]$. It follows from Definition \ref{def3.2} that $G_M(F(c_k,z))=(x\otimes y)\otimes z\neq x\otimes(y\otimes z)=G_M(F(x,F(y,z)))$. On the other hand, by associativity of $F$, we have $G_M(F(x,F(y,z)))=G_M(F(F(x,y),z))$. Thus $G_M(F(c_k,z))\neq G_M(F(F(x,y),z))$. Therefore, by Lemma \ref{lem3.3}, $(\min\{F(c_k,z),F(F(x,y),z)\},\max\{F(c_k,z),F(F(x,y),z)\})\cap (M\setminus\mathcal{C})\neq\emptyset $. Obviously, $I_{k}^{y}=O(\{c_k\}\cup F(M_{k}^{y},y))= \bigcup_{m,n\in \{c_k\}\cup F(M_{k}^{y},y)}(\min\{m,n\}, \max\{m,n\})$, $z\in M^{y}$, and $x\in M_{k}^{y}$. So that $(\min\{F(c_k,z),F(F(x,y),z)\},\max\{F(c_k,z),F(F(x,y),z)\})\subseteq F(I_{k}^{y},M^{y})$, which implies $F(I_{k}^{y},M^{y})\cap (M\setminus\mathcal{C})\neq\emptyset$.

(ii) The case that $F(x,y)\in M\setminus\mathcal{C}$ and $F(y,z)\notin M\setminus\mathcal{C}$ is completely analogous.

(iii) Let $F(x,y)\notin M\setminus\mathcal{C}$ and $F(y,z)\notin M\setminus\mathcal{C}$. Then from $F(x,y)\notin M\setminus\mathcal{C}$, there exists a $k\in K$ such that $F(x,y)\in [b_k,d_k]$. Thus $x\otimes y=c_k$ where $\{c_k\}=M\cap [b_k,d_k]$. Meanwhile, from $F(y,z)\notin M\setminus\mathcal{C}$, there exists an $l\in K$ such that $F(y,z)\in [b_l,d_l]$. Hence $y\otimes z=c_l$ where $\{c_l\}=M\cap [b_l,d_l]$. Consequently, $G_M(F(c_k,z))=(x\otimes y)\otimes z\neq x\otimes(y\otimes z)=G_M(F(x,c_l))$. Furthermore, by applying Lemma \ref{lem3.3}, we have $(\min\{F(c_k,z),F(x,c_l)\}, \max\{F(c_k,z),F(x,c_l)\})\cap (M\setminus\mathcal{C})\neq \emptyset$. Obviously, $x\in M_{k}^{y}$ and $F(x,c_l)\in F(M_{k}^{y},c_l)$. Similarly, $z\in M_{l}^{y}$ and $F(c_k,z)\in F(c_k,M_{l}^{y})$. Therefore, $$(\min\{F(c_k,z),F(x,c_l)\},\max\{F(c_k,z),F(x,c_l)\})\subseteq J_{k,l}^{y}.$$ This follows $J_{k,l}^{y}\cap (M\setminus\mathcal{C})\neq \emptyset$.

Cases (i), (ii) and (iii) deduce that $\mathfrak{T}(M)\cap (M\setminus\mathcal{C})\neq\emptyset$.

Conversely, suppose $\mathfrak{T}(M)\cap (M\setminus\mathcal{C})\neq\emptyset$. Then there exist a $y\in M$ and two elements $k,l\in K$ such that $F(I_{k}^{y},M^{y})\cap (M\setminus\mathcal{C})\neq\emptyset$ or $J_{k,l}^{y}\cap (M\setminus\mathcal{C})\neq \emptyset$. We distinguish two cases as follows.

(i) If $F(I_{k}^{y},M^{y})\cap (M\setminus\mathcal{C})\neq\emptyset$, then there exists a $z\in M^{y}$ such that $F(I_{k}^{y},z)\cap (M\setminus\mathcal{C})\neq\emptyset$. Thus by the definition of $I_{k}^{y}$, $F(M_{k}^{y},y)\neq\emptyset$. Applying Lemma \ref{lem3.2}, there exist elements $u\in\{c_k\}$ and $v\in F(M_{k}^{y},y)$ such that
$$(\min\{F(u,z),F(v,z)\},\max\{F(u,z),F(v,z)\})\cap (M\setminus\mathcal{C})\neq\emptyset.$$
Because of $v\in F(M_{k}^{y},y)$, there exists an $x\in M_{k}^{y}$ such that $F(x,y)=v$. Therefore, there exist elements $u\in\{c_k\}$ and $x\in M_{k}^{y}$ such that
$$(\min\{F(c_k,z),F(F(x,y),z)\},\max\{F(c_k,z),F(F(x,y),z)\})\cap (M\setminus\mathcal{C})\neq\emptyset.$$
Consequently, from Lemma \ref{lem3.3} we have $G_M(F(c_k,z))\neq G_M(F(F(x,y),z))$. On the other hand, from $x\in M_{k}^{y}$ we have $F(x,y)\in [b_k,d_k]$. Thus $x\otimes y=c_k$ where $\{c_k\}=M\cap [b_k,d_k]$. Moreover, from $z\in M_{y}$, we have $F(y,z)\in M\setminus\mathcal{C}$, and this follows $y\otimes z=F(y,z)$. Therefore, $(x\otimes y)\otimes z=G_M(F(c_k,z))\neq G_M(F(F(x,y),z))=G_M(F(x,F(y,z)))= x\otimes(y\otimes z)$.

(ii) If $J_{k,l}^{y}\cap (M\setminus\mathcal{C})\neq \emptyset$, then $J_{k,l}^{y}\neq \emptyset$. Thus by the definition of $J_{k,l}^{y}$, $F(O(F(M_{k}^{y},c_l)\cup F(c_k,M_{l}^{y})),1)\cap (M\setminus\mathcal{C})\neq \emptyset$ where $1$ is a neutral element of t-norm $F$, $F(M_{k}^{y},c_l)\neq\emptyset$ and $F(c_k,M_{y}^{l})\neq\emptyset$. Applying Lemma \ref{lem3.2}, there exist elements $u\in F(c_k,M_{l}^{y})$ and $v\in F(M_{k}^{y},c_l)$ such that $$(\min\{F(u,1),F(v,1)\},\max\{F(u,1),F(v,1)\})\cap (M\setminus\mathcal{C})\neq\emptyset.$$
Because $u\in F(c_k,M_{l}^{y})$ and $v\in F(M_{k}^{y},c_l)$, there exist an $x\in M_{k}^{y}$ and a $z\in M_{l}^{y}$ such that $u=F(c_k,z)$, $v=F(x,c_l)$. Therefore,
$$(\min\{F(c_k,z),F(x,c_l)\},\max\{F(c_k,z),F(x,c_l)\})\cap (M\setminus\mathcal{C})\neq\emptyset.$$
Hence, by Lemma \ref{lem3.3} we have $G_M(F(c_k,z))\neq G_M(F(x,c_l))$. On the other hand, from $x\in M_{k}^{y}$ we have $F(x,y)\in [b_k,d_k]$. Thus $x\otimes y=c_k$ where $\{c_k\}=M\cap [b_k,d_k]$.
Moreover, from $z\in M_{l}^{y}$ we have $F(y,z)\in [b_l,d_l]$, and this means $y\otimes z=c_l$ where $\{c_l\}=M\cap [b_l,d_l]$. Therefore, $(x\otimes y)\otimes z=G_M(F(c_k,z))\neq G_M(F(x,c_l))= x\otimes(y\otimes z)$.
\end{proof}

\begin{corollary}\label{coro3.1}
Let $M\in \mathcal{A}$ and $\mathfrak{T}(M)$ be the $F$-transformation of $M$. Then $(M,\otimes)$ is a semigroup if and only if $\mathfrak{T}(M)\cap ACC_{0}(M)=\emptyset$.
\end{corollary}

Furthermore, Proposition \ref{prop3.3} and Corollary \ref{coro3.1} imply the following theorem.
\begin{theorem}\label{thm3.2}
Let $f:[0,1]\rightarrow [0,1]$ be a strictly increasing function and $T:[0,1]^2\rightarrow [0,1]$ be a function defined by Eq.(\ref{eq1}). Then the function $T$ is associative if and only if $\mathfrak{T}(\emph{Ran}(f))\cap ACC_{0}(\emph{Ran}(f))=\emptyset$.
\end{theorem}

The following examples illustrate the above-mentioned results.
\begin{example}\label{exap3.1}
\renewcommand{\labelenumi}{(\roman{enumi})}
\begin{enumerate}
\item \emph{Let $F(x,y)=x\cdot y$ for all $x,y\in[0,1]$ and the function $f:[0,1]\rightarrow [0,1]$ be defined by
\begin{equation*}
f(x)=e^{x-1}.
\end{equation*}
Then $M=[e^{-1},1] $. It is clear that $\mathfrak{T}(M)\cap (M\setminus\mathcal{C})=\emptyset$. Hence, $T$ given by Eq.(\ref{eq1}) is an associative function by Theorem \ref{thm3.2}. Indeed, from Eq.(\ref{eq1}),
\begin{equation*}
T(x,y)=\max\{0,x+y-1\}
\end{equation*}
for all $x,y\in[0,1]$.
It is obvious that $T$ is a t-norm.}

\item \emph{Consider the function $F:[0,1]^2\rightarrow [0,1]$ given by \begin{equation*}
F(x,y)=\left\{
  \begin{array}{ll}
    0.5\cdot \max\{0,2x+2y-1\} & \hbox{if }(x,y)\in[0,0.5)^{2}, \\
     0.5+0.5\cdot (2x-1)(2y-1) & \hbox{if }(x,y)\in[0.5,1]^{2}, \\
    \min\{x,y\} & \hbox{otherwise}
  \end{array}
\right.
\end{equation*}
 and the strictly increasing function $f:[0,1]\rightarrow [0,1]$ defined by \begin{equation*}
 f(x)=\begin{cases}
h(x) & \hbox{if } x\leq0.5,\\
x &  \hbox{otherwise}
\end{cases}
\end{equation*}
where $h$ is a strictly increasing function with $h(0.5)\leq0.25$.
By a simple calculation, we get $M=\mbox{Ran}(h)\cup (0.5,1]$. It is a matter of straightforward verification to show that $\mathfrak{T}(M) \cap (M\setminus\mathcal{C})=\emptyset$. Therefore, $T$ given by Eq.(\ref{eq1}) is an associative function by Theorem \ref{thm3.2}. In fact, from Eq.(\ref{eq1}),
\begin{equation*}
T(x,y)=\left\{
  \begin{array}{ll}
    0 & \hbox{if }(x,y)\in[0,0.5]^{2}, \\
    0.5+0.5\cdot (2x-1)(2y-1) & \hbox{if }(x,y)\in(0.5,1]^{2}, \\
    \min\{x,y\} & \hbox{otherwise.}
  \end{array}
\right.
\end{equation*}
It is easy to check that $T$ is a t-norm.}

\item  \emph{In (ii), let $h(0.5)=0.5$. Then clearly $\mathfrak{T}(M) \cap (M\setminus\mathcal{C})=\emptyset$. Therefore, $T$ given by Eq.(\ref{eq1}) is an associative function by Theorem \ref{thm3.2}. In fact, from Eq.(\ref{eq1}),
\begin{equation*}
T(x,y)=\left\{
  \begin{array}{ll}
    0 & \hbox{if }(x,y)\in[0,0.5)^{2}, \\
    0.5+0.5\cdot (2x-1)(2y-1) & \hbox{if }(x,y)\in[0.5,1]^{2}, \\
    \min\{x,y\} & \hbox{otherwise.}
  \end{array}
\right.
\end{equation*}
Obviously, $T$ is a t-norm.}
\item \emph{Consider the function $F$ shown in (ii) again
 and the strictly increasing function $f:[0,1]\rightarrow [0,1]$ defined by
 \begin{equation*}
 f(x)=\begin{cases}
h(x) & \hbox{if } x\leq0.5,\\
0.5+0.5\cdot e^{2x-2} &  \hbox{otherwise}
\end{cases}
\end{equation*}
where $h$ is a strictly increasing function with $h(0.5)\leq0.25$.
 It is easy to check that  $\mathfrak{T}(M) \cap (M\setminus\mathcal{C})\neq\emptyset$ and $T$ given by Eq.(\ref{eq1}) is not an associative function by Theorem \ref{thm3.2}. In fact, from Eq.(\ref{eq1}),
\begin{equation*}
T(x,y)=\left\{
  \begin{array}{ll}
    0 & \hbox{if }(x,y)\in[0,0.5]^{2}, \\
    0.5+0.5\cdot \max\{0,2x+2y-3\} & \hbox{if }(x,y)\in(0.5,1]^{2}, \\
    \min\{x,y\} & \hbox{otherwise.}
  \end{array}
\right.
\end{equation*}
Taking $x=y=0.75,z=0.5$, we have $T(x,y)=T(0.75,0.75)=0.5$, $T(x,z)=T(y,z)=T(0.75,0.5)=0.5$ and $T(0.5,0.5)=0$. Hence $T(T(x,y),z)=0\neq 0.5=T(x,T(y,z))$.}
\end{enumerate}
\end{example}

As shown by Examples \ref{exap3.1} (ii),(iii) and (iv), for a given t-norm $F$ that is an ordinal sum of t-norms in the sense of A. H. Clifford, the function $T$ obtained by \eqref{eq1} may not be necessarily associative when we slightly modify the strictly increasing function $f$. Next, we further establish the conditions under which the function $T$ given by \eqref{eq1} is associative when the t-norm $F$ is an ordinal sum of t-norms in the sense of A. H. Clifford.

\section{Characterizations of the associativity of function $T$ given by \eqref{eq1} when the t-norm $F$ is an ordinal sum of t-norms}

In this section, we characterize the associativity of function $T$ given by \eqref{eq1} when the t-norm $F$ is an ordinal sum of t-norms in the sense of A. H. Clifford.

Let $A\neq\emptyset$ be a totally ordered set and $F=(\langle a_\alpha,b_\alpha, F_\alpha\rangle)_{\alpha\in A}$ be an ordinal sum of a family of t-norms $\{F_\alpha \mid \alpha\in A\}$. Then for each $\alpha\in A$, the t-norms $F_\alpha:[0,1]^2\rightarrow [0,1]$ and $F^\alpha:[a_\alpha,b_\alpha]^2\rightarrow [a_\alpha,b_\alpha]$ can be constructed by $$F_\alpha(x,y)=\frac{F((a_\alpha+(b_\alpha-a_\alpha)x),(a_\alpha+(b_\alpha-a_\alpha)y))-a_\alpha}{b_\alpha-a_\alpha}$$ and
$$F^\alpha(x,y)=a_\alpha+(b_\alpha-a_\alpha)F_\alpha(\frac{x-a_\alpha}{b_\alpha-a_\alpha},\frac{y-a_\alpha}{b_\alpha-a_\alpha}),$$ respectively. Note that if the t-norm $F$ has an idempotent element $x\in(a_\alpha,b_\alpha)$ with $\alpha\in A$ then there is a $y\in(a_\alpha,b_\alpha)$ such that $F(x,y)<y$. Moreover, let $f:[0,1]\rightarrow [0,1]$ be a strictly increasing function with $\mbox{Ran}(f)\subseteq[0,1]$ and denote \begin{equation*}
s_\alpha=\inf\{x\in [0,1]\mid f(x)\geq a_\alpha\mid \alpha\in A\}
\end{equation*}
and
\begin{equation*}
t_\alpha=\sup\{x\in [0,1]\mid f(x)\leq b_\alpha\mid \alpha\in A\}.
\end{equation*}
Thus, if there exists $\alpha\in A$ such that $\mbox{Ran}(f)\cap[a_\alpha, b_\alpha]$ is an empty set or a singleton set then $s_\alpha=t_\alpha$ and $(s_\alpha,t_\alpha)=\emptyset$. If there exists $\alpha\in A$ such that $\mbox{Ran}(f)\cap[a_\alpha, b_\alpha]$ is infinite then we have two results as follows.
\begin{enumerate}
\item If there exists an $x\in[0,1]$ such that $f(x)=a_\alpha$, then $s_\alpha=x$. Otherwise, there exists a $y\in[0,1]$ such that $a_\alpha\leq f(y^+)< b_\alpha$ and either that $f(y)>a_\alpha$ and $f(y^-)\leq a_\alpha$ (we just have $f(0)>a_\alpha$ if $y=0$) or that $f(y)<a_\alpha$, in this case, $s_\alpha=y$.
\item If there exists an $x\in[0,1]$ such that $f(x)=b_\alpha$, then $t_\alpha=x$. Otherwise, there exists a $y\in[0,1]$ such that $a_\alpha<f(y^-)\leq b_\alpha$ and either that $f(y)<b_\alpha$ and $f(y^+)\geq b_\alpha$ (we just have $f(1)<b_\alpha$ if $y=1$) or that $f(y)>b_\alpha$, in this case, $t_\alpha=y$.
\end{enumerate}
In particular, if $\mbox{Ran}(f)\cap[a_\alpha, b_\alpha]$ with $\alpha\in A$ is an empty set or a singleton set, then the function $T$ given by \eqref{eq1} coincides with $T_M(x,y)=\min\{x,y\}$ for all $x,y\in[0,1]$. Therefore, in what follows, we always suppose that there exists an $\alpha\in A$ such that $\mbox{Ran}(f)\cap[a_\alpha, b_\alpha]$ is infinite.

By collecting all the non-empty intervals $(s_\alpha,t_\alpha)_{\alpha\in A}$ and their corresponding intervals $(a_\alpha,b_\alpha)_{\alpha\in A}$, we define a new index set $B_A^f\neq\emptyset$ with $B_A^f\subseteq A$ and the order inherited from $A$.
 %and consider the function $o:B_A^f\rightarrow[0,1]$, which assigns the midpoint to each of the intervals under consideration by $$o(\beta)=\frac{s_\beta+t_\beta}{2}$$ where $\beta\in B_A^f$. It is easy to see that $o$ is injective and the intervals $(s_\alpha,t_\alpha)_{\alpha\in A}$ cannot be overlapping.
 Obviously, $\{(s_\beta,t_\beta)\}_{\beta\in B_A^f}$ is a family of pairwise disjoint open subintervals of $[0,1]$. For every $\beta\in B_A^f$, define a function $f_\beta:[s_\beta,t_\beta]\rightarrow [a_\beta,b_\beta]$ by
\begin{equation}\label{eq06}
f_\beta(x)=\left\{
  \begin{array}{ll}
    f(x) & \hbox{if } x\in(s_\beta,t_\beta), \\
    f(s_\beta) & \hbox{if } x=s_\beta \hbox{ and } f(s_\beta)\geq a_\beta, \\
    a_\beta & \hbox{if } x=s_\beta \hbox{ and } f(s_\beta)< a_\beta, \\
    f(t_\beta) & \hbox{if } x=t_\beta \hbox{ and } f(t_\beta)\leq b_\beta, \\
    b_\beta & \hbox{if } x=t_\beta \hbox{ and } f(t_\beta)> b_\beta.
  \end{array}
\right.
\end{equation}
Then $f_\beta$ is a strictly increasing function. The set $\{f_\beta\mid f_\beta:[s_\beta,t_\beta]\rightarrow [a_\beta,b_\beta]\}_{\beta\in B_A^f}$ is called a \emph{decomposition set} of $f$.

We first have the following proposition.
\begin{proposition}\label{prop4.1}
Let $A\neq\emptyset$ be a totally ordered set, $F=(\langle a_\alpha,b_\alpha, F_\alpha\rangle)_{\alpha\in A}$ be an ordinal sum of a family of t-norms $\{F_\alpha \mid \alpha\in A\}$, $f$ be a strictly increasing function and the set $\{f_\beta\mid f_\beta:[s_\beta,t_\beta]\rightarrow [a_\beta,b_\beta]\}_{\beta\in B_A^f}$ be a decomposition set of $f$. Then $T(x,y)=T^\beta(x,y)$ for all $x,y\in[s_{\beta},t_{\beta}]$ and $\beta\in B_A^f$ except $x=y=s_\beta$ where $T(x,y)=f^{(-1)}(F(f(x),f(y)))$ and $T^\beta(x,y)=f_\beta^{(-1)}(F^\beta(f_\beta(x),f_\beta(y)))$.
\end{proposition}
\begin{proof}Obviously, $f^{(-1)}(y)=f_\beta^{(-1)}(y)$ for all $y\in[a_\beta,b_\beta]$ and $f_\beta(x)=f(x)$ for all $x\in(s_\beta,t_\beta)$ from the definition of function $f_\beta$. Then $T(x,y)=f^{(-1)}(F(f(x),f(y)))=f_\beta^{(-1)}(F^\beta(f_\beta(x),f_\beta(y)))=T^\beta(x,y)$ for all $x,y\in(s_\beta,t_\beta)$ and $\beta\in B_A^f$. The remaining proof is split into three cases.
\renewcommand{\labelenumi}{(\roman{enumi})}
\begin{enumerate}
\item The case when $f(s_\beta)\geq a_\beta$ or $f(t_\beta)\leq b_\beta$. If $f(s_\beta)\geq a_\beta$ then $T(s_\beta,y)=T^\beta(s_\beta,y)$ for all $y\in(s_\beta,t_\beta]$ is obvious. If $f(t_\beta)\leq b_\beta$, then $T(t_\beta,y)=T^\beta(t_\beta,y)$ for all $y\in[s_\beta,t_\beta]$ is obvious.

\item The case when $f(s_\beta)<a_\beta$.
For any $y\in(s_\beta,t_\beta]$, we have
 $$T(s_\beta,y)=f^{(-1)}(F(f(s_\beta),f(y)))=f^{(-1)}(\min \{f(s_\beta),f(y)\})=f^{(-1)}(f(s_\beta))=s_\beta$$
 since $f(s_\beta)<a_\beta<f(y)$.
 On the other hand, from the definition of function $f_\beta$ we have
$$T^\beta(s_\beta,y)=f_\beta^{(-1)}(F^\beta(f_\beta(s_\beta),f_\beta(y)))=f_\beta^{(-1)}(F^\beta(a_\beta,f_\beta(y)))=f_\beta^{(-1)}(a_\beta)=f_\beta^{(-1)}(f_\beta(s_\beta))=s_\beta$$ for all $y\in[s_\beta,t_\beta]$.
Therefore, $T(s_\beta,y)=T^\beta(s_\beta,y)$ for all $y\in(s_\beta,t_\beta]$.
\item The case when $f(t_\beta)>b_\beta$. We get for all $y\in[s_\beta,t_\beta)$, $$T(t_\beta,y)=f^{(-1)}(F(f(t_\beta),f(y)))=f^{(-1)}(\min\{f(t_\beta),f(y)\})=f^{(-1)}(f(y))=y.$$ On the other hand, $F(f(t_\beta),f(t_\beta))\leq f(t_\beta)$ since $F$ is a t-norm, and $F(f(t_\beta),f(t_\beta))\geq b_\beta\geq f(t_\beta^{-})$. Thus $T(t_\beta,t_\beta)=f^{(-1)}(F(f(t_\beta),f(t_\beta)))=t_\beta$. We also get that for all $y\in[s_\beta,t_\beta]$, $$T^\beta(t_\beta,y)=f_\beta^{(-1)}(F^\beta(f_\beta(t_\beta),f_\beta(y)))=f_\beta^{(-1)}(F^\beta(b_\beta,f_\beta(y)))=f_\beta^{(-1)}(f_\beta(y))=y.$$ Therefore, $T(t_\beta,y)=T^\beta(t_\beta,y)$ for all $y\in[s_\beta,t_\beta]$.
\end{enumerate}
The commutativity of $F$ implies that for any $\beta\in B_A^f$, $T(x,y)=T^\beta(x,y)$ for all $x,y\in[s_{\beta},t_{\beta}]$ except $x=y=s_\beta$.
\end{proof}

Then we have the following proposition.
\begin{proposition}\label{prop4.2}
Let $A\neq\emptyset$ be a totally ordered set, $F=(\langle a_\alpha,b_\alpha, F_\alpha\rangle)_{\alpha\in A}$ be an ordinal sum of a family of t-norms $\{F_\alpha \mid \alpha\in A\}$, $f$ be a strictly increasing function and the set $\{f_\beta\mid f_\beta:[s_\beta,t_\beta]\rightarrow [a_\beta,b_\beta]\}_{\beta\in B_A^f}$ be a decomposition set of $f$. Then $T^\beta$ is associative on $[s_\beta,t_\beta]^2$  if and only if $\mathfrak{T}_\beta(\emph{Ran}(f_\beta))\cap ACC_{0}(\emph{Ran}(f_\beta))=\emptyset$ where the set $\mathfrak{T}_\beta(\emph{Ran}(f_\beta))$ is the $F^\beta$-transformation of $\emph{ Ran}(f_\beta)$ and $T^\beta(x,y)=f_\beta^{(-1)}(F^\beta(f_\beta(x),f_\beta(y)))$.
\end{proposition}
\begin{proof} By Theorem \ref{thm3.2}, $T^\beta(x,y)$ is associative on $[s_\beta,t_\beta]^2$ if and only if $\mathfrak{T}_\beta(\mbox{Ran}(f_\beta))\cap ACC_{0}(\mbox{Ran}(f_\beta))=\emptyset$.
\end{proof}

Now we introduce the following theorem cited from \cite{SJ2002,EP2002}.
\begin{theorem}\label{theorem:4.1}\emph{(Ordinal Sum Theorem for t-subnorms)}.
Let $A\neq\emptyset$ be a totally ordered set, $(G_\alpha)_{\alpha\in A}$ be a family of t-subnorms and $\{(a_\alpha,b_\alpha)\}_{\alpha\in A}$ be a family of pairwise disjoint open subintervals of $[0, 1]$. Further, if $b_{\alpha_0}=a_{\alpha_1}$ for some $\alpha_0,\alpha_1\in A$ then assume either that $G_{\alpha_0}$ is a t-norm or that $G_{\alpha_1}$ has no zero divisors. Let $G:[0,1]^2\rightarrow[0,1]$ be a function defined by \begin{equation}\label{eq:4.1}
G(x,y)=\left\{
  \begin{array}{ll}
    a_\alpha+(b_\alpha-a_\alpha)G_\alpha(\frac{x-a_\alpha}{b_\alpha-a_\alpha},\frac{y-a_\alpha}{b_\alpha-a_\alpha}) & \hbox{if }(x,y)\in(a_\alpha,b_\alpha]^{2}, \\
    \min\{x,y\} & \hbox{otherwise.}
  \end{array}
\right.
\end{equation}
Then $G$ is a t-subnorm that is known as an ordinal sum of $\{([a_\alpha,b_\alpha],G_\alpha) \mid \alpha\in A\}$, and is written as $G=(\langle a_\alpha,b_\alpha, G_\alpha\rangle)_{\alpha\in A}$. Further, $G_\alpha$and $(a_\alpha,b_\alpha)$ are referred to as a \emph{summand} and an \emph{open generating interval} of the ordinal sum, respectively.
\end{theorem}

Furthermore, if $b_{\alpha_2}=1$ for some $\alpha_2\in A$ then assume that $G_{\alpha_2}$ is a t-norm. Then the $G$ defined by Eq.\eqref{eq:4.1} is a t-norm. The construction in Theorem \ref{theorem:4.1} is the most general way to obtain a t-subnorm (resp. a t-norm) as an ordinal sum of semigroups.

\begin{proposition}\label{prop4.3}
Let $p, q, u, v\in [-\infty, \infty]$ with $p<q, u<v$ and a function $T: [p,q]^2\rightarrow [p,q]$ be defined by $T(x,y)=f^{(-1)}(F(f(x),f(y)))$ where $f: [p,q]\rightarrow [u,v]$ is a strictly increasing function and $F:[u,v]^2\rightarrow[u,v]$ is a t-norm.
Then there exists a $t\in(p,q]$ such that $T(t,x)=x$ if and only if $f(x^-)\leq F(f(t),f(x))$ for all $x\leq t$.
\end{proposition}
\begin{proof}
Suppose there exists a $t\in(p,q]$ such that $T(t,x)=x$ for all $x\leq t$. Then $f^{(-1)}(F(f(t),f(x)))=T(t,x)=x=f^{(-1)}(f(x))$. Thus, from $F(f(t),f(x))\leq f(x)$, $f^{(-1)}(F(f(t),f(x)))=f^{(-1)}(f(x))$ if and only if $\mbox{Ran}(f)\cap[F(f(t),f(x)),f(x)]$ contains at most one element if and only if $f(x^-)\leq F(f(t),f(x))$. Therefore, $T(t,x)=x$ if and only if $f(x^-)\leq F(f(t),f(x))$.
\end{proof}

It is easy to see that the function $T$ defined by \eqref{eq1} is always commutative, non-decreasing and $T(x, y)\leq \min\{x, y\}$ for all $x, y\in[0, 1]$. Further, we have the following proposition.
\begin{proposition}\label{prop4.4}
Let $A\neq\emptyset$ be a totally ordered set, $F=(\langle a_\alpha,b_\alpha, F_\alpha\rangle)_{\alpha\in A}$ be an ordinal sum of a family of t-norms $\{F_\alpha \mid \alpha\in A\}$, $f$ be a strictly increasing function and the set $\{f_\beta\mid f_\beta:[s_\beta,t_\beta]\rightarrow [a_\beta,b_\beta]\}_{\beta\in B_A^f}$ be a decomposition set of $f$. If the following two statements are hold:
\renewcommand{\labelenumi}{(\roman{enumi})}
\begin{enumerate}
\item $\mathfrak{T}_\beta(\emph{Ran}(f_\beta))\cap ACC_{0}(\emph{Ran}(f_\beta))=\emptyset$ for any $\beta\in B_A^f$, and
\item if $t_{\beta_i}=s_{\beta_j}$ for some $\beta_i,\beta_j\in B_A^f$ and there are two elements $u,v\in (s_{\beta_j},t_{\beta_j})$ such that $$F(f_{\beta_j}(u),f_{\beta_j}(v))\leq f(s_{\beta_j}^+),$$ then $f(x^-)\leq F(f(t_{\beta_i}),f(x))$ for all $x\leq t_{\beta_i}$,
\end{enumerate}
then the function $T$ defined by \eqref{eq1} is a t-subnorm.
\end{proposition}
\begin{proof} First note that by Proposition \ref{prop4.1}, $T(x,y)=T^\beta(x,y)$ for all $x,y\in[s_{\beta},t_{\beta}]$ and $\beta\in B_A^f$ except $x=y=s_\beta$.
Now, for each $\beta\in B_A^f$, construct $G_\beta:[0,1]^2\rightarrow [0,1]$ and $G^\beta:[0,1]^2\rightarrow [0,1]$ by
\begin{equation*}
G_\beta(x,y)=\left\{
  \begin{array}{ll}
   \frac{T((s_{\beta}+(t_{\beta}-s_{\beta})x),(s_{\beta}+(t_{\beta}-s_{\beta})y))-s_{\beta}}{t_{\beta}-s_{\beta}} & \hbox{either }T(s_\beta,s_\beta)=s_\beta\\ &\hbox{ or } \min\{x,y\}\neq 0 \hbox{ and }T(s_\beta,s_\beta)<s_\beta,\\
   0 & \hbox{if } \min\{x,y\}=0 \hbox{ and }T(s_\beta,s_\beta)<s_\beta
  \end{array}
\right.
\end{equation*}
and
\begin{equation*}
G^\beta(x,y)=\frac{T^\beta((s_{\beta}+(t_{\beta}-s_{\beta})x),(s_{\beta}+(t_{\beta}-s_{\beta})y))-s_{\beta}}{t_{\beta}-s_{\beta}},
\end{equation*}
respectively.
Then one can check $G_\beta=G^\beta$. This together with the statement (i) and Proposition \ref{prop4.2} implies that $G_\beta$ is a t-subnorm.

 On the other hand, if $t_{\beta_i}=s_{\beta_j}$ for some $\beta_i,\beta_j\in B_A^f$ and there exist two elements $u,v\in (s_{\beta_j},t_{\beta_j})$ such that $F(f_{\beta_j}(u),f_{\beta_j}(v))\leq f(s_{\beta_j}^+)$, then we get that $T(u,v)=s_{\beta_j}$, i.e., $G_{\beta_j}$ has zero divisors. Furthermore, by the statement (ii) we have $f(x^-)\leq F(f(t_{\beta_i}),f(x))$ for all $x\leq t_{\beta_i}$, which together Proposition \ref{prop4.3} deduces that $T(t_{\beta_i},x)=x$ for all $x\leq t_{\beta_i}$. Hence, $t_{\beta_i}$ is the neutral element of $T|_{[s_{\beta_i},t_{\beta_i}]^2}$ where $T|_{[s_{\beta_i},t_{\beta_i}]^2}$ is a restriction of the function $T$ on $[s_{\beta_i},t_{\beta_i}]^2$, i.e., $G_{\beta_i}$ is a t-norm.

 Therefore, there is a family of t-subnorms $\{G_\beta\}_{\beta\in B_A^f}$ such that either $G_{\beta_i}$ is a t-norm or $G_{\beta_j}$ has no zero divisors when $t_{\beta_i}=s_{\beta_j}$ for some $\beta_i,\beta_j\in B_A^f$, and
$$T(x,y)=s_{\beta}+(t_{\beta}-s_{\beta})G_\beta(\frac{x-s_{\beta}}{t_{\beta}-s_{\beta}},\frac{y-s_{\beta}}{t_{\beta}-s_{\beta}})\mbox{ for all }(x,y)\in(s_\beta,t_\beta]^{2}.$$

In order to complete the proof, we just need to examine the structure of function $T$ on $[0,1]^2\setminus \bigcup_{\beta\in B_A^f}(s_{\beta},t_{\beta}]^2$. There are five cases as follows.
\renewcommand{\labelenumi}{(\roman{enumi})}
\begin{enumerate}
\item If $x\in[0,1]\setminus \bigcup_{\beta\in B_A^f}[s_{\beta},t_{\beta}]$ and $y\in[0,1]$, then we have $$T(x,y)=f^{(-1)}(F(f(x),f(y)))=f^{(-1)}(\min\{f(x),f(y)\})=\min\{x,y\}.$$

\item If there exists a $\beta\in B_A^f$ such that $x\in[s_{\beta},t_{\beta}]$, then there are three subcases as follows.
\begin{enumerate}
\item If $x\in(s_{\beta},t_{\beta})$, then $T(x,y)=f^{(-1)}(F(f(x),f(y)))=f^{(-1)}(\min\{f(x),f(y)\})=\min\{x,y\}$ for all $y\in[0,1]\setminus [s_{\beta},t_{\beta}]$.
\item If $x=s_{\beta}$, then there are two subcases as below.
\begin{itemize}
\item If $[s_{\beta},t_{\beta}]\cap\bigcup_{\{\gamma\in B\mid \gamma<\beta\}}[s_{\gamma},t_{\gamma}]=\emptyset$, then $$T(s_{\beta},y)=f^{(-1)}(F(f(s_{\beta}),f(y)))=f^{(-1)}(\min\{f(s_{\beta}),f(y)\})=\min\{s_{\beta},y\}$$ for all $y\in[0,1]\setminus [s_{\beta},t_{\beta}]$.
\item If there exists a $\gamma_0\in B_A^f$ such that $[s_{\beta},t_{\beta}]\cap\bigcup_{\{\gamma\in B\mid \gamma<\beta\}}[s_{\gamma},t_{\gamma}]=\{s_{\beta}\}=\{t_{\gamma_0}\}$, then $T(s_{\beta},y)=f^{(-1)}(F(f(s_{\beta}),f(y)))=f^{(-1)}(\min\{f(s_{\beta}),f(y)\})=\min\{s_{\beta},y\}$ for all $y\in[0,1]\setminus [s_{\gamma_0},t_{\beta}]$.
\end{itemize}
\item If $x=t_{\beta}$, then there are two subcases as follows.
\begin{itemize}
\item If $[s_{\beta},t_{\beta}]\cap\bigcup_{\{\gamma\in B\mid \gamma>\beta\}}[s_{\gamma},t_{\gamma}]=\emptyset$, then $$T(t_{\beta},y)=f^{(-1)}(F(f(t_{\beta}),f(y)))=f^{(-1)}(\min\{f(t_{\beta}),f(y)\})=\min\{t_{\beta},y\}$$ for all $y\in[0,1]\setminus [s_{\beta},t_{\beta}]$.
\item If there exists a $\gamma_1\in B_A^f$ such that $[s_{\beta},t_{\beta}]\cap\bigcup_{\{\gamma\in B\mid \gamma>\beta\}}[s_{\gamma},t_{\gamma}]=\{t_{\beta}\}=\{s_{\gamma_1}\}$, then $T(t_{\beta},y)=f^{(-1)}(F(f(t_{\beta}),f(y)))=f^{(-1)}(\min\{f(t_{\beta}),f(y)\})=\min\{t_{\beta},y\}$ for all $y\in[0,1]\setminus [s_{\beta},t_{\gamma_1}]$.
\end{itemize}
\end{enumerate}

\item If there is a $\beta\in B_A^f$ such that $x=s_{\beta}$ and $y\in(s_{\beta},t_{\beta}]$, then we distinguish two subcases:
\begin{enumerate}
\item If $f(s_\beta)\geq a_\beta$, then $f(s_{\beta}^-)\leq a_{\beta}\leq F(f(s_{\beta}),f(y))\leq f(s_{\beta})$. Hence, $$T(s_{\beta},y)=f^{(-1)}(F(f(s_{\beta}),f(y)))=s_{\beta}=\min\{s_{\beta},y\}.$$
\item If $f(s_\beta)<a_\beta$, then $$T(s_\beta,y)=f^{(-1)}(F(f(s_\beta),f(y)))=f^{(-1)}(\min \{f(s_\beta),f(y)\})=f^{(-1)}(f(s_\beta))=s_\beta=\min\{s_{\beta},y\}$$ since $f(y)> a_\beta$.
\end{enumerate}
Therefore, in any case, $T(s_{\beta},y)=\min\{s_{\beta},y\}$ for any $y\in(s_{\beta},t_{\beta}]$.
\item If there is a $\beta\in B_A^f$ such that $x\in(s_{\beta},t_{\beta}]$ and $y=s_{\beta}$, then in complete analogy to (iii), we have $T(x,s_{\beta})=\min\{x,s_{\beta}\}$.
\item If $x=y=s_\beta$, then there are two subcases as follows.
\renewcommand{\labelenumi}{(\roman{enumi})}
\begin{enumerate}
\item If $[s_{\beta},t_{\beta}]\cap\bigcup_{\{\gamma\in B\mid \gamma<\beta\}}[s_{\gamma},t_{\gamma}]=\emptyset$, then we have two subcases as below.
\renewcommand{\labelenumi}{(\roman{enumi})}
\begin{enumerate}
\item If $f(s_\beta)\geq a_\beta$, then $T(s_{\beta},s_{\beta})=f^{(-1)}(F(f(s_{\beta}),f(s_{\beta})))=s_{\beta}=\min\{s_{\beta},s_{\beta}\}$.
\item If $f(s_\beta)< a_\beta$, then there are two subcases as follows.
\begin{itemize}
\item If $\{f(s_\beta)\}\cap \bigcup_{\alpha\in A}(a_\alpha,b_\alpha)=\emptyset$, then $f(s_\beta)$ is an idempotent element of $F$. Thus $T(s_\beta,s_\beta)=f^{(-1)}(F(f(s_\beta),f(s_\beta)))=f^{(-1)}(f(s_\beta))=s_\beta=\min\{s_\beta,s_{\beta}\}$.
\item  If there exists an $\alpha_0\in A$ such that $\{f(s_\beta)\}\cap \bigcup_{\alpha\in A}(a_\alpha,b_\alpha)=\{f(s_\beta)\}$, then $\mbox{Ran}(f)\cap[a_{\alpha_0},b_{\alpha_0}]$ is a singleton set and we get $f(s_\beta^-)\leq a_{\alpha_0} \leq F(f(s_\beta),f(s_\beta))\leq f(s_\beta)$. Therefore, $T(s_\beta,s_\beta)=f^{(-1)}(F(f(s_\beta),f(s_\beta)))=s_\beta=\min\{s_\beta,s_{\beta}\}$.
\end{itemize}
\end{enumerate}
\item If there exists a $\gamma_0\in B_A^f$ such that $[s_{\beta},t_{\beta}]\cap\bigcup_{\{\gamma\in B\mid \gamma<\beta\}}[s_{\gamma},t_{\gamma}]=\{s_{\beta}\}=\{t_{\gamma_0}\}$, then
$T(s_{\beta},s_{\beta})=T(t_{\gamma_0},t_{\gamma_0})\in \{T(x,y)\mid x,y\in(s_{\beta_i},t_{\beta_i}]\}_{\beta_i\in B}$.
\end{enumerate}
\end{enumerate}

Consequently,
\begin{equation*}
T(x,y)=\left\{
  \begin{array}{ll}
    s_{\beta}+(t_{\beta}-s_{\beta})G_\beta(\frac{x-s_{\beta}}{t_{\beta}-s_{\beta}},\frac{y-s_{\beta}}{t_{\beta}-s_{\beta}}) & \hbox{if }(x,y)\in(s_\beta,t_\beta]^{2}, \\
    \min\{x,y\} & \hbox{otherwise.}
  \end{array}
\right.
\end{equation*}

Finally, from Theorem \ref{theorem:4.1} the function $T$ defined by \eqref{eq1} is an ordinal sum of a countable number of semigroups and has the associativity.
\end{proof}

\begin{proposition}\label{prop4.5}
Let $A\neq\emptyset$ be a totally ordered set, $F=(\langle a_\alpha,b_\alpha, F_\alpha\rangle)_{\alpha\in A}$ be an ordinal sum of a family of t-norms $\{F_\alpha \mid \alpha\in A\}$, $f$ be a strictly increasing function and the set $\{f_\beta\mid f_\beta:[s_\beta,t_\beta]\rightarrow [a_\beta,b_\beta]\}_{\beta\in B_A^f}$ be a decomposition set of $f$. Then $\mathfrak{T}(\emph{Ran}(f))\cap ACC_{0}(\emph{Ran}(f))=\emptyset$ if and only if the following two statements are hold:
\renewcommand{\labelenumi}{(\roman{enumi})}
\begin{enumerate}
\item $\mathfrak{T}_\beta(\emph{Ran}(f_\beta))\cap ACC_{0}(\emph{Ran}(f_\beta))=\emptyset$ for any $\beta\in B_A^f$;
\item if $t_{\beta_i}=s_{\beta_j}$ for some $\beta_i,\beta_j\in B_A^f$ and there are two elements $u,v\in (s_{\beta_j},t_{\beta_j})$ such that $$F(f_{\beta_j}(u),f_{\beta_j}(v))\leq f(s_{\beta_j}^+),$$ then $f(x^-)\leq F(f(t_{\beta_i}),f(x))$ for all $x\leq t_{\beta_i}$.
\end{enumerate}
\end{proposition}
\begin{proof} If the statements (i) and (ii) hold, then from Theorem \ref{thm3.2} and Proposition \ref{prop4.4}, we have $\mathfrak{T}(\mbox{Ran}(f))\cap ACC_{0}(\mbox{Ran}(f))=\emptyset$.

Conversely, assume that $\mathfrak{T}(\mbox{Ran}(f))\cap ACC_{0}(\mbox{Ran}(f))=\emptyset$ holds. Then Propositions \ref{prop4.1} and \ref{prop4.2} mean that (i) is true. If $t_{\beta_i}=s_{\beta_j}$ for some $\beta_i,\beta_j\in B_A^f$ and there exist two elements $u,v\in (s_{\beta_j},t_{\beta_j})$ such that $F(f_{\beta_j}(u),f_{\beta_j}(v))\leq f(s_{\beta_j}^+)$, then by Theorem \ref{thm3.2}, $T(t_{\beta_i},z)=T(T(u,v),z)=T(u,T(v,z))=T(u,\min\{v,z\})=T(u,z)=\min\{u,z\}=z$ for any $z< t_{\beta_i}$. This together with Proposition \ref{prop4.3} results in $f(z^-)\leq F(f(t_{\beta_i}),f(z))$ for all $z< t_{\beta_i}$. Now, let $z=t_{\beta_i}$. Then $T(t_{\beta_i},t_{\beta_i})=T(T(u,v),z)=T(u,T(v,z))=T(u,T(v,t_{\beta_i}))=T(u,T(v,s_{\beta_j}))=T(u,s_{\beta_j})=s_{\beta_j}=t_{\beta_i}$, i.e.,  $T(t_{\beta_i},t_{\beta_i})=t_{\beta_i}$. This follows that $f(t_{\beta_i}^-)\leq F(f(t_{\beta_i}),f(t_{\beta_i}))$. Therefore, if $t_{\beta_i}=s_{\beta_j}$ for some $\beta_i,\beta_j\in B_A^f$ and there exist $u,v\in (s_{\beta_j},t_{\beta_j})$ such that $F(f_{\beta_j}(u),f_{\beta_j}(v))\leq f(s_{\beta_j}^+)$ then $f(x^-)\leq F(f(t_{\beta_i}),f(x))$ for all $x\leq t_{\beta_i}$.
\end{proof}

\begin{remark}\label{rem:4.2} \emph{If $t_{\beta_i}=s_{\beta_j}$ for some $\beta_i,\beta_j\in B_A^f$ and there exist elements $u,v\in (s_{\beta_j},t_{\beta_j})$ and $p\in(s_{\beta_i}, t_{\beta_i}]$ such that $F(f_{\beta_j}(u),f_{\beta_j}(v))\leq f(s_{\beta_j}^+)$ and $F(f(t_{\beta_i}),f(p))<f(p^-)$, then $\mathfrak{T_{1}}(M) \cap (M\setminus\mathcal{C})\neq\emptyset$ or $\mathfrak{T_{2}}(M) \cap (M\setminus\mathcal{C})\neq\emptyset$.}

\emph{Indeed, let $f_{\beta_j}(u)=x$, $f_{\beta_j}(v)=y$ and $f_{\beta_i}(p)=z$.
Then there are two cases as follows.
    \begin{enumerate}
\item If $p$ is a continuous point of $f$, then $F(x,y)\notin M\setminus\mathcal{C}$ and $F(y,z)\in M\setminus\mathcal{C}$. Hence, $y\otimes z=F(y,z)$ and there exists a $k\in K$ such that $F(x,y)\in [b_k,d_k]$. Thus $x\otimes y=c_k$ where $\{c_k\}=M\cap[b_k,d_k]$. On the other hand, from Proposition \ref{prop4.5} and Corollary \ref{coro3.1}, we know that $(M,\otimes)$ is not a semigroup. Thus it follows from Definition \ref{def3.2} that $G_M(F(c_k,z))=(x\otimes y)\otimes z\neq x\otimes(y\otimes z)=G_M(F(x,F(y,z)))$. Moreover, the associativity of $F$ implies $G_M(F(x,F(y,z)))=G_M(F(F(x,y),z))$. Therefore, $G_M(F(c_k,z))\neq G_M(F(F(x,y),z))$. Applying Lemma \ref{lem3.3}, we have $$(\min\{F(c_k,z),F(F(x,y),z)\},\max\{F(c_k,z),F(F(x,y),z)\})\cap (M\setminus\mathcal{C})\neq\emptyset.$$ So that $(\min\{F(c_k,z),F(F(x,y),z)\},\max\{F(c_k,z),F(F(x,y),z)\})\subseteq F(I_{k}^{y},M^{y})$, which implies $F(I_{k}^{y},M^{y})\cap (M\setminus\mathcal{C})\neq\emptyset$, i.e., $\mathfrak{T_{1}}(M) \cap (M\setminus\mathcal{C})\neq\emptyset$.
\item If $p$ is a discontinuous point of $f$, then $F(x,y)\notin M\setminus\mathcal{C}$ and $F(y,z)\notin M\setminus\mathcal{C}$. Hence from $F(x,y)\notin M\setminus\mathcal{C}$, there exists a $k\in K$ such that $F(x,y)\in [b_k,d_k]$, thus $x\otimes y=c_k$ where $\{c_k\}=M\cap [b_k,d_k]$. From $F(y,z)\notin M\setminus\mathcal{C}$, there exists an $l\in K$ such that $F(y,z)\in [b_l,d_l]$, then $y\otimes z=c_l$ where $\{c_l\}=M\cap [b_l,d_l]$. On the other hand, from Proposition \ref{prop4.5} and Corollary \ref{coro3.1} we know that $(M,\otimes)$ is not a semigroup. Therefore, $G_M(F(c_k,z))=(x\otimes y)\otimes z\neq x\otimes(y\otimes z)=G_M(F(x,c_l))$. Applying Lemma \ref{lem3.3}, we have $$(\min\{F(c_k,z),F(x,c_l)\}, \max\{F(c_k,z),F(x,c_l)\})\cap (M\setminus\mathcal{C})\neq \emptyset.$$ Obviously, $x\in M_{k}^{y}$ and $F(x,c_l)\in F(M_{k}^{y},c_l)$. Similarly, $z\in M_{l}^{y}$ and $F(c_k,z)\in F(c_k,M_{l}^{y})$. Therefore, $$(\min\{F(c_k,z),F(x,c_l)\},\max\{F(c_k,z),F(x,c_l)\})\subseteq J_{k,l}^{y}.$$ This follows $J_{k,l}^{y}\cap (M\setminus\mathcal{C})\neq \emptyset$, i.e., $\mathfrak{T_{2}}(M) \cap (M\setminus\mathcal{C})\neq\emptyset$.
\end{enumerate}}
\end{remark}

\begin{theorem}\label{thm4.2}
Let $A\neq\emptyset$ be a totally ordered set, $F=(\langle a_\alpha,b_\alpha, F_\alpha\rangle)_{\alpha\in A}$ be an ordinal sum of a family of t-norms $\{F_\alpha \mid \alpha\in A\}$, $f$ be a strictly increasing function and the set $\{f_\beta\mid f_\beta:[s_\beta,t_\beta]\rightarrow [a_\beta,b_\beta]\}_{\beta\in B_A^f}$ be a decomposition set of $f$. Then the function $T:[0,1]^2\rightarrow[0,1]$ defined by \eqref{eq1} is a t-norm if and only if the following three statements are satisfied:
\renewcommand{\labelenumi}{(\roman{enumi})}
\begin{enumerate}
\item $\mathfrak{T}_\beta(\emph{Ran}(f_\beta))\cap ACC_{0}(\emph{Ran}(f_\beta))=\emptyset$ for any $\beta\in B_A^f$;
\item if $t_{\beta_i}=s_{\beta_j}$ for some $\beta_i,\beta_j\in B_A^f$ and there exist two elements $u,v\in (s_{\beta_j},t_{\beta_j})$ such that $$F(f_{\beta_j}(u),f_{\beta_j}(v))\leq f(s_{\beta_j}^+),$$ then $f(x^-)\leq F(f(t_{\beta_i}),f(x))$ for all $x\leq t_{\beta_i}$;
\item $f(x^-)\leq F(f(1),f(x))$ for all $x\leq 1$.
\end{enumerate}
\end{theorem}
\begin{proof} By Proposition \ref{prop4.5}, the statements (i) and (ii) are true if and only if the function $T$ given by Eq.\eqref{eq1} is a t-subnorm. From Proposition \ref{prop4.3}, $f(x^-)\leq F(f(1),f(x))$ for all $x\leq 1$ if and only if $T(1,x)=x$ for all $x\leq 1$. Therefore, the function $T$ is a t-norm if and only if 1 is a neutral element of $T$.
\end{proof}

\begin{example}\label{exap4.1}
\renewcommand{\labelenumi}{(\roman{enumi})}
\begin{enumerate}
\item \emph{ In Example \ref{exap3.1} (iv), $[s_{\beta_1},t_{\beta_1}]=[0,0.5]$, $[s_{\beta_2},t_{\beta_2}]=[0.5,1]$. From Eq.\eqref{eq06}, $f_{\beta_1}:[0,0.5]\rightarrow [0,0.5]$ is defined by $f_{\beta_1}(x)=h(x)$ and $f_{\beta_2}:[0.5,1]\rightarrow [0.5,1]$ is defined by $f_{\beta_2}(x)=0.5+0.5\cdot e^{2x-2}$. Clearly, there exist two elements $u,v\in (0.5,1)$ such that $F(f_{\beta_2}(u),f_{\beta_2}(v))\leq f(0.5^+)$. However, $ F(f(0.5),f(x))=0< f(x^-)$ for all $0<x\leq0.5$. Therefore, the function $T$ defined by Eq.\eqref{eq1} is not associative by Proposition \ref{prop4.5}.}

\item \emph{Consider the function $F$ given by \begin{equation*}
F(x,y)=\left\{
  \begin{array}{ll}
    0.5\cdot T^{nM}(2x,2y)& \hbox{if }(x,y)\in[0,0.5)^{2}, \\
     0.5+0.5\cdot (2x-1)(2y-1) & \hbox{if }(x,y)\in[0.5,1]^{2}, \\
    \min\{x,y\} & \hbox{otherwise}
  \end{array}
\right.
\end{equation*}
where $T^{nM}:[0,1]^2\rightarrow[0,1]$ is defined by \begin{equation*}
T^{nM}(x,y)=\left\{
  \begin{array}{ll}
    0& \hbox{if } x+y\leq 1, \\
     \min\{x,y\} &\hbox{otherwise}
  \end{array}
\right.
\end{equation*}
 and the strictly increasing function $f:[0,1]\rightarrow [0,1]$ defined by \begin{equation*}
 f(x)=\begin{cases}
0.2x+0.3 & \hbox{if } x\in[0,0.5],\\
0.5+0.5\cdot e^{2x-2} &  \hbox{otherwise}.
\end{cases}
\end{equation*}
Then $[s_{\beta_1},t_{\beta_1}]=[0,0.5]$, $[s_{\beta_2},t_{\beta_2}]=[0.5,1]$. It is clear that $\mathfrak{T}_{\beta_i}(\mbox{Ran}(f_{\beta_i}))\cap ACC_{0}(\mbox{Ran}(f_{\beta_i}))=\emptyset$ for $i=1,2$. For any $0<x\leq 0.5$, from $2f(x)+2f(0.5)>1$ we have $F(f(x),f(0.5))=0.5T^{nM}(2f(x),2f(0.5))=0.5\cdot\min\{2f(x),2f(0.5)\}=\min\{f(x),f(0.5)\}=f(x)\geq f(x^-)$. Then by Proposition \ref{prop4.5}, the function $T$ defined by Eq.(\ref{eq1}) is associative. In fact, from Eq.(\ref{eq1}),
\begin{equation*}
T(x,y)=\left\{
  \begin{array}{ll}
    0.5+0.5\cdot \max\{2x+2y-3,0\} & \hbox{if }(x,y)\in(0.5,1]^{2}, \\
    \min\{x,y\} & \hbox{otherwise.}
  \end{array}
\right.
\end{equation*}
It is easy to check that $T$ is a t-norm.}
\item \emph{Consider the function $F$ given by \begin{equation*}
F(x,y)=\left\{
  \begin{array}{ll}
    0.5\cdot (2x)(2y) & \hbox{if }(x,y)\in[0,0.5]^{2}, \\
     0.5+0.5\cdot (2x-1)(2y-1) & \hbox{if }(x,y)\in[0.5,1]^{2}, \\
    \min\{x,y\} & \hbox{otherwise}
  \end{array}
\right.
\end{equation*}
 and the strictly increasing function $f:[0,1]\rightarrow [0,1]$ defined by \begin{equation*}
 f(x)=\begin{cases}
0.8x & \hbox{if } x<0.5,\\
0.5+0.5\cdot e^{2x-2} &  \hbox{otherwise}.
\end{cases}
\end{equation*}
Then $[s_{\beta_1},t_{\beta_1}]=[0,0.5]$, $[s_{\beta_2},t_{\beta_2}]=[0.5,1]$. It is clear that $\mathfrak{T}_\beta(\mbox{Ran}(f_{\beta_i}))\cap ACC_{0}(\mbox{Ran}(f_{\beta_i}))=\emptyset$ for $i=1,2$. For any $x< 0.5$, from $f(x)<0.5<f(0.5)$ we have $F(f(x),f(0.5))=\min\{f(x),f(0.5)\}=f(x)\geq f(x^-)$. Then by Proposition \ref{prop4.5}, the function $T$ defined by Eq.(\ref{eq1}) is associative. In fact, from Eq.(\ref{eq1}),
\begin{equation*}
T(x,y)=\left\{
  \begin{array}{ll}
    \frac{8}{5}xy & \hbox{if }(x,y)\in[0,0.5)^{2}, \\
    0.5+0.5\cdot \max\{2x+2y-3,0\} & \hbox{if }(x,y)\in(0.5,1]^{2}, \\
    \min\{x,y\} & \hbox{otherwise.}
  \end{array}
\right.
\end{equation*}}
\end{enumerate}
\end{example}

\section{Characterizations of the associativity of function $T$ given by \eqref{eq1} when the t-norm $F$ is an ordinal sum of t-subnorms}

In this section, we characterize the associativity of function $T$ given by \eqref{eq1} when the t-norm $F$ is an ordinal sum of t-subnorms in the sense of A. H. Clifford.

In the proof of Theorem \ref{them3.1}, we used the fact that 1 is a neutral element of t-norm $F$. In general, a t-subnorm has not a neutral element. Therefore, Theorem \ref{them3.1} is not suitable for a proper t-subnorm.

In the following, we temporarily turn to give a necessary and sufficient condition for the function $T$ given by \eqref{eq1} being associative when $F$ is a proper t-subnorm.

Let $F:[0,1]^2\rightarrow[0,1]$ be a proper t-subnorm and $f: [0,1]\rightarrow [0,1]$ be a strictly increasing function. Define $\overline{F}:[0,1]^2\rightarrow[0,1]$ by \begin{equation}\label{eq5.1}
\overline{F}(x,y)=\left\{
  \begin{array}{ll}
    \frac{1}{2}\cdot F(2x,2y) & \hbox{if }(x,y)\in[0,0.5]^{2}, \\
    \min\{x,y\} & \hbox{otherwise}
  \end{array}
\right.
\end{equation}
and $\overline{f}: [0,1]\rightarrow [0,1]$ by
$$\overline{f}(x)=\frac{1}{2}f(x)$$
for all $x\in[0,1]$, respectively.
Then $\overline{F}:[0,1]^2\rightarrow[0,1]$ is obviously a t-norm as an ordinal sum of t-suborms and $\overline{f}:[0,1]\rightarrow[0,1]$ is a strictly increasing function.

\begin{proposition}\label{prop5.1}
Let $F:[0,1]^2\rightarrow[0,1]$ be a proper t-subnorm, $\overline{F}$ be define as Eq. \eqref{eq5.1}, $f: [0,1]\rightarrow [0,1]$ be a strictly increasing function and $\overline{f}(x)=\frac{1}{2}f(x)$ for all $x\in[0,1]$. Then $T(x,y)=\overline{T}(x,y)$ for all $x,y\in[0,1]$ where $T(x,y)=f^{(-1)}(F(f(x),f(y)))$ and $\overline{T}(x,y)=\overline{f}^{(-1)}(\overline{F}(\overline{f}(x),\overline{f}(y)))$.
\end{proposition}
\begin{proof} Since $\overline{f}(x)=\frac{1}{2}f(x)$ for all $x\in[0,1]$, we have $f^{(-1)}(y)=\overline{f}^{(-1)}(\frac{y}{2})$ for all $y\in[0,1]$. Therefore, for all $x,y\in[0,1]$,
\begin{eqnarray*}
\overline{T}(x,y)&=&\overline{f}^{(-1)}(\overline{F}(\overline{f}(x),\overline{f}(y)))\\
&=&\overline{f}^{(-1)}(\frac{1}{2}F(2\overline{f}(x),2\overline{f}(y)))\\
&=&\overline{f}^{(-1)}(\frac{1}{2}F(f(x),f(y))\\
&=&f^{(-1)}(F(f(x),f(y)))\\
&=&T(x,y).
\end{eqnarray*}
\end{proof}

In particular, Theorem \ref{thm3.2} and Proposition \ref{prop5.1} imply the following theorem.
\begin{theorem}\label{them5.1}
Let $F:[0,1]^2\rightarrow[0,1]$ be a proper t-subnorm, $\overline{F}$ be define as Eq. \eqref{eq5.1}, $f: [0,1]\rightarrow [0,1]$ be a strictly increasing function and $\overline{f}(x)=\frac{1}{2}f(x)$ for all $x\in[0,1]$. Then the function $T$ given by \eqref{eq1} is associative if and only if $\mathfrak{\overline{T}}(\emph{Ran}(\overline{f}))\cap ACC_{0}(\emph{Ran}(\overline{f}))=\emptyset$, where $\mathfrak{\overline{T}}(\emph{Ran}(\overline{f}))$ is the $\overline{F}$-transformation of $\emph{Ran}(\overline{f})$.
\end{theorem}

Next, we come back to discuss the associativity of function $T$ given by \eqref{eq1} when the t-norm $F$ is an ordinal sum of t-subnorms.

Let $A\neq\emptyset$ be a totally ordered set, and the t-norm $F=(\langle a_\alpha,b_\alpha, F_\alpha\rangle)_{\alpha\in A}$ be an ordinal sum of a family of t-subnorms $\{F_\alpha \mid \alpha\in A\}$.
In the following, denote \begin{equation*}\label{eq5.2}
F^\alpha(x,y)=a_\alpha+(b_\alpha-a_\alpha)F_\alpha(\frac{x-a_\alpha}{b_\alpha-a_\alpha},\frac{y-a_\alpha}{b_\alpha-a_\alpha})
\end{equation*}
for all $x,y\in[a_\alpha,b_\alpha]$.
Then the only possible difference between $F^\alpha$ and $F$ on the interval $[a_\alpha,b_\alpha]^2$ is the value of the point at $(a_\alpha,a_\alpha)$. In fact, if $b_{\alpha_0}=a_{\alpha_1}$ for some $\alpha_0,\alpha_1\in A$ and $F_{\alpha_1}$ has zero divisors, then $F(a_{\alpha_1},a_{\alpha_1})=a_{\alpha_1}=F^\alpha(a_{\alpha_1},a_{\alpha_1})$ and $F|_{[a_{\alpha_1},b_{\alpha_1}]^2}$ is a t-subnorm.
If $b_{\alpha_0}=a_{\alpha_1}$ for some $\alpha_0,\alpha_1\in A$ and $F_{\alpha_1}$ has no zero divisors, then it is possible that $F(a_{\alpha_1},a_{\alpha_1})=F(b_{\alpha_0},b_{\alpha_0})<b_{\alpha_0}= a_{\alpha_1}=F^{\alpha_1}(a_{\alpha_1},a_{\alpha_1})$.

We have the following proposition.
\begin{proposition}\label{prop5.2}
Let $A\neq\emptyset$ be a totally ordered set, the t-norm $F=(\langle a_\alpha,b_\alpha, F_\alpha\rangle)_{\alpha\in A}$ be an ordinal sum of a family of t-subnorms $\{F_\alpha \mid \alpha\in A\}$, $f$ be a strictly increasing function and the set $\{f_\beta\mid f_\beta:[s_\beta,t_\beta]\rightarrow [a_\beta,b_\beta]\}_{\beta\in B_A^f}$ be a decomposition set of $f$.
\renewcommand{\labelenumi}{(\roman{enumi})}
\begin{enumerate}
\item If $F^\beta$ is a t-norm on $[a_{\beta},b_{\beta}]^2$, then $T(x,y)=T^\beta(x,y)$ for all $x,y\in[s_{\beta},t_{\beta}]$ except $x=y=s_\beta$, where $T(x,y)=f^{(-1)}(F(f(x),f(y)))$ and $T^\beta(x,y)=f_\beta^{(-1)}(F^\beta(f_\beta(x),f_\beta(y)))$.
\item If $F^\beta$ is a proper t-subnorm on $[a_{\beta},b_{\beta}]^2$, then
\renewcommand{\labelenumi}{(\roman{enumi})}
\begin{enumerate}
\item if $f(t_{\beta})\leq b_{\beta}$ then $T(x,y)=T^\beta(x,y)$ for all $x,y\in[s_{\beta},t_{\beta}]$ except $x=y=s_\beta$.
\item if $f(t_{\beta})> b_{\beta}$ then $T(x,y)=\underline{T}^\beta(x,y)$ for all $x,y\in[s_{\beta},t_{\beta}]$ except $x=y=s_\beta$,
where $\underline{T}^\beta(x,y)=f_\beta^{(-1)}(\underline{F}^\beta(f_\beta(x),f_\beta(y)))$ and \begin{equation*}
\underline{F}^\beta(x,y)=\left\{
  \begin{array}{ll}
    F^\beta(x,y) & \emph{if }(x,y)\in[a_{\beta},b_{\beta})^{2}, \\
    \min\{x,y\} & \emph{if }(x,y)\in[a_{\beta},b_{\beta}]^{2}\setminus[a_{\beta},b_{\beta})^{2}.
  \end{array}
\right.
\end{equation*}
\end{enumerate}
\end{enumerate}
\end{proposition}
\begin{proof} (i) It is immediately by Proposition \ref{prop4.1}.

(ii) If $F^\beta$ is a proper t-subnorm on $[a_{\beta},b_{\beta}]^2$, then (a) is obvious. So that we just prove (b).

In fact, if $f(t_\beta)>b_\beta$, then by the definition of function $f_\beta$ we have $f_\beta(t_\beta)=b_\beta$ and $$T(t_\beta,y)=f^{(-1)}(F(f(t_\beta),f(y)))=f^{(-1)}(\min\{f(t_\beta),f(y)\})=f^{(-1)}(f(y))=y$$ for all $y\in[s_\beta,t_\beta)$. In particular, if $y=t_\beta$, then $F(f(t_\beta),f(y))=F(f(t_\beta),f(t_\beta))\leq f(t_\beta)$ since $F$ is a t-norm. The last inequality together with $F(f(t_\beta),f(t_\beta))\geq b_\beta\geq f(t_\beta^{-})$ implies $T(t_\beta,t_\beta)=f^{(-1)}(F(f(t_\beta),f(t_\beta)))=t_\beta$. Thus $T(t_\beta,y)=y$ for all $y\in[s_\beta,t_\beta]$. Moreover, $$\underline{T}^\beta(t_\beta,y)=f_\beta^{(-1)}(\underline{F}^\beta(f_\beta(t_\beta),f_\beta(y)))=f_\beta^{(-1)}(\underline{F}^\beta(b_\beta,f_\beta(y)))=f_\beta^{(-1)}(f_\beta(y))=y$$ for all $y\in[s_\beta,t_\beta]$. Therefore, $T(t_\beta,y)=\underline{T}^\beta(t_\beta,y)$ for all $y\in[s_\beta,t_\beta]$.

 The case $(x,y)\in [s_\beta,t_\beta)^2\setminus \{(s_\beta,s_\beta)\}$ is obvious. Hence, $T(x,y)=\underline{T}^\beta(x,y)$ for all $x,y\in[s_\beta,t_\beta]$ except $x=y=s_\beta$.
\end{proof}

\begin{proposition}\label{prop5.3}
Let $A\neq\emptyset$ be a totally ordered set, the t-norm $F=(\langle a_\alpha,b_\alpha, F_\alpha\rangle)_{\alpha\in A}$ be an ordinal sum of a family of t-subnorms $\{F_\alpha \mid \alpha\in A\}$, $f$ be a strictly increasing function and the set $\{f_\beta\mid f_\beta:[s_\beta,t_\beta]\rightarrow [a_\beta,b_\beta]\}_{\beta\in B_A^f}$ be a decomposition set of $f$.
\renewcommand{\labelenumi}{(\roman{enumi})}
\begin{enumerate}
\item If $F^\beta$ is a t-norm on $[a_{\beta},b_{\beta}]^2$, then
$$T^\beta \emph{ is associative on } [s_\beta,t_\beta]^2 \emph{ if and only if } \mathfrak{T}_\beta(\emph{Ran}(f_\beta))\cap ACC_{0}(\emph{Ran}(f_\beta))=\emptyset,$$
where $T^\beta(x,y)=f_\beta^{(-1)}(F^\beta(f_\beta(x),f_\beta(y)))$ and $\mathfrak{T}_\beta(\emph{Ran}(f_\beta))$ is the $F^\beta$-transformation of $\emph{ Ran}(f_\beta)$.
\item If $F^\beta$ is a proper t-subnorm on $[a_{\beta},b_{\beta}]^2$, then \renewcommand{\labelenumi}{(\roman{enumi})}
\begin{enumerate}
\item if $f(t_{\beta})\leq b_{\beta}$ then $$
T^\beta \emph{ is associative on } [s_\beta,t_\beta]^2 \emph{ if and only if } \mathfrak{\overline{T}}_\beta(\emph{Ran}(\overline{f}_\beta))\cap ACC_{0}(\emph{Ran}(\overline{f}_\beta))=\emptyset,
$$
where $\mathfrak{\overline{T}}_\beta(\emph{Ran}(\overline{f}_\beta))$ is the $\overline{F}^\beta$-transformation of $\emph{ Ran}(\overline{f}_\beta)$,
$\overline{f}_\beta(x)=\frac{a_\beta+f_\beta(x)}{2}$ for any $x\in[s_\beta,t_\beta]$ and
\begin{equation*}
\overline{F}^\beta(x,y)=\left\{
  \begin{array}{ll}
    a_\beta+\frac{(b_\beta-a_\beta)}{2}F_\beta(\frac{2(x-a_\beta)}{b_\beta-a_\beta},\frac{2(y-a_\beta)}{b_\beta-a_\beta}) & \emph{if }(x,y)\in[a_{\beta},\frac{b_{\beta}+a_{\beta}}{2}]^{2}, \\
    \min\{x,y\} & \emph{if }(x,y)\in[a_{\beta},b_{\beta}]^{2}\setminus[a_{\beta},\frac{b_{\beta}+a_{\beta}}{2}]^{2}.
  \end{array}
\right.
\end{equation*}

\item if $f(t_{\beta})> b_{\beta}$ then $$
\underline{T}^\beta \emph{ is associative on } [s_\beta,t_\beta]^2 \emph{ if and only if } \mathfrak{\underline{T}}_\beta(\emph{Ran}(f_\beta))\cap ACC_{0}(\emph{Ran}(f_\beta))=\emptyset,$$
where $\underline{T}^\beta(x,y)=f_\beta^{(-1)}(\underline{F}^\beta(f_\beta(x),f_\beta(y)))$, $\mathfrak{\underline{T}}_\beta(\emph{Ran}(f_\beta))$ is the $\underline{F}^\beta$-transformation of $\emph{ Ran}(f_\beta)$ and \begin{equation*}
\underline{F}^\beta(x,y)=\left\{
  \begin{array}{ll}
    F^\beta(x,y) & \emph{if }(x,y)\in[a_{\beta},b_{\beta})^{2}, \\
    \min\{x,y\} & \emph{if }(x,y)\in[a_{\beta},b_{\beta}]^{2}\setminus[a_{\beta},b_{\beta})^{2}.
  \end{array}
\right.
\end{equation*}
\end{enumerate}
\end{enumerate}
\end{proposition}
\begin{proof} (i) It is immediate from Proposition \ref{prop4.2}.

(ii)(a) If $F^\beta$ is a proper t-subnorm on $[a_{\beta},b_{\beta}]^2$ and $f(t_{\beta})\leq b_{\beta}$, then let $\overline{f}_\beta(x)=\frac{a_\beta+f_\beta(x)}{2}$ for any $x\in[s_\beta,t_\beta]$ and \begin{equation*}
\overline{F}^\beta(x,y)=\left\{
  \begin{array}{ll}
    a_\beta+\frac{(b_\beta-a_\beta)}{2}F_\beta(\frac{2(x-a_\beta)}{b_\beta-a_\beta},\frac{2(y-a_\beta)}{b_\beta-a_\beta}) & \hbox{if }(x,y)\in[a_{\beta},\frac{b_{\beta}+a_{\beta}}{2}]^{2}, \\
    \min\{x,y\} & \hbox{if }(x,y)\in[a_{\beta},b_{\beta}]^{2}\setminus[a_{\beta},\frac{b_{\beta}+a_{\beta}}{2}]^{2}.
  \end{array}
\right.
\end{equation*} Thus one can check that $f_\beta(x)=2\overline{f}_\beta(x)-a_\beta$ for all $x\in[s_\beta,t_\beta]$ and $f_\beta^{(-1)}(y)=\overline{f}_\beta^{(-1)}(\frac{y+a_\beta}{2})$ for all $y\in[a_\beta,b_\beta]$. So that
 \begin{eqnarray*}
\overline{T}^\beta(x,y)&=&\overline{f}_\beta^{(-1)}(\overline{F}^\beta(\overline{f}_\beta(x),\overline{f}_\beta(y)))\\
&=&\overline{f}_\beta^{(-1)}(a_\beta+\frac{b_\beta-a_\beta}{2}F_\beta(\frac{2(\overline{f}_\beta(x)-a_\beta)}{b_\beta-a_\beta},\frac{2(\overline{f}_\beta(y)-a_\beta)}{b_\beta-a_\beta}))\\
&=&\overline{f}_\beta^{(-1)}(a_\beta+\frac{b_\beta-a_\beta}{2}F_\beta(\frac{f_\beta(x)-a_\beta}{b_\beta-a_\beta},\frac{f_\beta(y)-a_\beta}{b_\beta-a_\beta}))\\
&=&f_\beta^{(-1)}(a_\beta+(b_\beta-a_\beta)F_\beta(\frac{f_\beta(x)-a_\beta}{b_\beta-a_\beta},\frac{f_\beta(y)-a_\beta}{b_\beta-a_\beta}))\\
&=&f_\beta^{(-1)}(F^\beta(f_\beta(x),f_\beta(y)))\\
&=&T^\beta(x,y)
\end{eqnarray*}
since $\overline{T}^\beta(x,y)=\overline{f}_\beta^{(-1)}(\overline{F}^\beta(\overline{f}_\beta(x),\overline{f}_\beta(y)))$ for any $x,y\in[s_\beta,t_\beta]$, in which $\overline{F}^\beta$ is a t-norm on $[a_{\beta},b_{\beta}]^2$ since $F^\beta$ is a proper t-subnorm on $[a_{\beta},b_{\beta}]^2$. Obviously, $\overline{f}_\beta$ is a strictly increasing function from $[s_\beta,t_\beta]$ to $[a_{\beta},b_{\beta}]$. Thus it follows from Theorem \ref{thm3.2} that $\overline{T}^\beta$ is associative on $[s_\beta,t_\beta]^2$ if and only if $\mathfrak{\overline{T}}_\beta(\mbox{Ran}(\overline{f}_\beta))\cap ACC_{0}(\mbox{Ran}(\overline{f}_\beta))=\emptyset$. Therefore, $T^\beta$ is associative on $[s_\beta,t_\beta]^2$ if and only if $$\mathfrak{\overline{T}}_\beta(\mbox{Ran}(\overline{f}_\beta))\cap ACC_{0}(\mbox{Ran}(\overline{f}_\beta))=\emptyset.$$

(ii)(b) It is clear that $\underline{F}^\beta$ is a t-norm on $[a_{\beta},b_{\beta}]$. Therefore, from Theorem \ref{thm3.2}, $\underline{T}^\beta$ is associative on $[s_\beta,t_\beta]^2$ if and only if $\mathfrak{\underline{T}}_\beta(\mbox{Ran}(f_\beta))\cap ACC_{0}(\mbox{Ran}(f_\beta))=\emptyset$.
\end{proof}

\begin{proposition}\label{prop5.4}
Let $A\neq\emptyset$ be a totally ordered set, the t-norm $F=(\langle a_\alpha,b_\alpha, F_\alpha\rangle)_{\alpha\in A}$ be an ordinal sum of a family of t-subnorms $\{F_\alpha \mid \alpha\in A\}$, $f$ be a strictly increasing function and the set $\{f_\beta\mid f_\beta:[s_\beta,t_\beta]\rightarrow [a_\beta,b_\beta]\}_{\beta\in B_A^f}$ be a decomposition set of $f$. If the following three statements hold:
\renewcommand{\labelenumi}{(\roman{enumi})}
\begin{enumerate}
\item if $F^\beta$ is a t-norm on $[a_{\beta},b_{\beta}]^2$, then $\mathfrak{T}_\beta(\emph{Ran}(f_\beta))\cap ACC_{0}(\emph{Ran}(f_\beta))=\emptyset$
where $\mathfrak{T}_\beta(\emph{Ran}(f_\beta))$ is the $F^\beta$-transformation of $\emph{ Ran}(f_\beta)$;
\item if $F^\beta$ is a proper t-subnorm on $[a_{\beta},b_{\beta}]^2$, then
\begin{enumerate}
\item
if $f(t_{\beta})\leq b_{\beta}$ then $ \mathfrak{\overline{T}}_\beta(\emph{Ran}(\overline{f}_\beta))\cap ACC_{0}(\emph{Ran}(\overline{f}_\beta))=\emptyset$
where the set $\mathfrak{\overline{T}}_\beta(\emph{Ran}(\overline{f}_\beta))$ is the $\overline{F}^\beta$-transformation of $\emph{ Ran}(\overline{f}_\beta)$,
$\overline{f}_\beta(x)=\frac{a_\beta+f_\beta(x)}{2}$ for any $x\in[s_\beta,t_\beta]$ and
\begin{equation*}
\overline{F}^\beta(x,y)=\left\{
  \begin{array}{ll}
    a_\beta+\frac{(b_\beta-a_\beta)}{2}F_\beta(\frac{2(x-a_\beta)}{b_\beta-a_\beta},\frac{2(y-a_\beta)}{b_\beta-a_\beta}) & \emph{if }(x,y)\in[a_{\beta},\frac{b_{\beta}+a_{\beta}}{2}]^{2}, \\
    \min\{x,y\} & \emph{if }(x,y)\in[a_{\beta},b_{\beta}]^{2}\setminus[a_{\beta},\frac{b_{\beta}+a_{\beta}}{2}]^{2};
  \end{array}
\right.
\end{equation*}

\item if $f(t_{\beta})> b_{\beta}$ then $\mathfrak{\underline{T}}_\beta(\emph{Ran}(f_\beta))\cap ACC_{0}(\emph{Ran}(f_\beta))=\emptyset$
where the set $\mathfrak{\underline{T}}_\beta(\emph{Ran}(f_\beta))$ is the $\underline{F}^\beta$-transformation of $\emph{ Ran}(f_\beta)$ and
\begin{equation*}
\underline{F}^\beta(x,y)=\left\{
  \begin{array}{ll}
    F^\beta(x,y) & \emph{if }(x,y)\in[a_{\beta},b_{\beta})^{2}, \\
    \min\{x,y\} & \emph{if }(x,y)\in[a_{\beta},b_{\beta}]^{2}\setminus[a_{\beta},b_{\beta})^{2};
  \end{array}
\right.
\end{equation*}
\end{enumerate}
\item if $t_{\beta_i}=s_{\beta_j}$ for some $\beta_i,\beta_j\in B_A^f$ and there exist two elements $u,v\in (s_{\beta_j},t_{\beta_j})$ such that $$F(f_{\beta_j}(u),f_{\beta_j}(v))\leq f(s_{\beta_j}^+),$$ then $f(x^-)\leq F(f(t_{\beta_i}),f(x))$ for all $x\leq t_{\beta_i}$;
\end{enumerate}
then the function $T:[0,1]^2\rightarrow[0,1]$ defined by \eqref{eq1} is a t-subnorm.
\end{proposition}
\begin{proof} In complete analogy to Proposition \ref{prop4.4}.
\end{proof}

\begin{proposition}\label{prop5.5}
Let $A\neq\emptyset$ be a totally ordered set, the t-norm $F=(\langle a_\alpha,b_\alpha, F_\alpha\rangle)_{\alpha\in A}$ be an ordinal sum of a family of t-subnorms $\{F_\alpha \mid \alpha\in A\}$, $f$ be a strictly increasing function and the set $\{f_\beta\mid f_\beta:[s_\beta,t_\beta]\rightarrow [a_\beta,b_\beta]\}_{\beta\in B_A^f}$ be a decomposition set of $f$. Then $\mathfrak{T}(\emph{Ran}(f))\cap ACC_{0}(\emph{Ran}(f))=\emptyset$ if and only if the following three statements hold:\renewcommand{\labelenumi}{(\roman{enumi})}
\begin{enumerate}
\item if $F^\beta$ is a t-norm on $[a_{\beta},b_{\beta}]^2$, then $\mathfrak{T}_\beta(\emph{Ran}(f_\beta))\cap ACC_{0}(\emph{Ran}(f_\beta))=\emptyset$
where $\mathfrak{T}_\beta(\emph{Ran}(f_\beta))$ is the $F^\beta$-transformation of $\emph{ Ran}(f_\beta)$;
\item if $F^\beta$ is a proper t-subnorm on $[a_{\beta},b_{\beta}]^2$, then
\begin{enumerate}
\item
if $f(t_{\beta})\leq b_{\beta}$ then $ \mathfrak{\overline{T}}_\beta(\emph{Ran}(\overline{f}_\beta))\cap ACC_{0}(\emph{Ran}(\overline{f}_\beta))=\emptyset$
where the set $\mathfrak{\overline{T}}_\beta(\emph{Ran}(\overline{f}_\beta))$ is the $\overline{F}^\beta$-transformation of $\emph{ Ran}(\overline{f}_\beta)$,
$\overline{f}_\beta(x)=\frac{a_\beta+f_\beta(x)}{2}$ for any $x\in[s_\beta,t_\beta]$ and
\begin{equation*}
\overline{F}^\beta(x,y)=\left\{
  \begin{array}{ll}
    a_\beta+\frac{(b_\beta-a_\beta)}{2}F_\beta(\frac{2(x-a_\beta)}{b_\beta-a_\beta},\frac{2(y-a_\beta)}{b_\beta-a_\beta}) & \emph{if }(x,y)\in[a_{\beta},\frac{b_{\beta}+a_{\beta}}{2}]^{2}, \\
    \min\{x,y\} & \emph{if }(x,y)\in[a_{\beta},b_{\beta}]^{2}\setminus[a_{\beta},\frac{b_{\beta}+a_{\beta}}{2}]^{2};
  \end{array}
\right.
\end{equation*}

\item if $f(t_{\beta})> b_{\beta}$ then $\mathfrak{\underline{T}}_\beta(\emph{Ran}(f_\beta))\cap ACC_{0}(\emph{Ran}(f_\beta))=\emptyset$
where the set $\mathfrak{\underline{T}}_\beta(\emph{Ran}(f_\beta))$ is the $\underline{F}^\beta$-transformation of $\emph{ Ran}(f_\beta)$ and
\begin{equation*}
\underline{F}^\beta(x,y)=\left\{
  \begin{array}{ll}
    F^\beta(x,y) & \emph{if }(x,y)\in[a_{\beta},b_{\beta})^{2}, \\
    \min\{x,y\} & \emph{if }(x,y)\in[a_{\beta},b_{\beta}]^{2}\setminus[a_{\beta},b_{\beta})^{2};
  \end{array}
\right.
\end{equation*}
\end{enumerate}
\item if $t_{\beta_i}=s_{\beta_j}$ for some $\beta_i,\beta_j\in B_A^f$ and there exist two elements $u,v\in (s_{\beta_j},t_{\beta_j})$ such that $$F(f_{\beta_j}(u),f_{\beta_j}(v))\leq f(s_{\beta_j}^+),$$ then $f(x^-)\leq F(f(t_{\beta_i}),f(x))$ for all $x\leq t_{\beta_i}$.
\end{enumerate}
\end{proposition}
\begin{proof} In complete analogy to Proposition \ref{prop4.5}.
\end{proof}

With Propositions \ref{prop4.3} and \ref{prop5.5} we have the following theorem.
\begin{theorem}\label{thm5.2}
Let $A\neq\emptyset$ be a totally ordered set, the t-norm $F=(\langle a_\alpha,b_\alpha, F_\alpha\rangle)_{\alpha\in A}$ be an ordinal sum of a family of t-subnorms $\{F_\alpha \mid \alpha\in A\}$, $f$ be a strictly increasing function and the set $\{f_\beta\mid f_\beta:[s_\beta,t_\beta]\rightarrow [a_\beta,b_\beta]\}_{\beta\in B_A^f}$ be a decomposition set of $f$. Then the function $T:[0,1]^2\rightarrow[0,1]$ defined by \eqref{eq1} is a t-norm if and only if the following four statements hold:\renewcommand{\labelenumi}{(\roman{enumi})}
\begin{enumerate}
\item if $F^\beta$ is a t-norm on $[a_{\beta},b_{\beta}]^2$, then $\mathfrak{T}_\beta(\emph{Ran}(f_\beta))\cap ACC_{0}(\emph{Ran}(f_\beta))=\emptyset$
where $\mathfrak{T}_\beta(\emph{Ran}(f_\beta))$ is the $F^\beta$-transformation of $\emph{ Ran}(f_\beta)$;
\item if $F^\beta$ is a proper t-subnorm on $[a_{\beta},b_{\beta}]^2$, then
\begin{enumerate}
\item
if $f(t_{\beta})\leq b_{\beta}$ then $ \mathfrak{\overline{T}}_\beta(\emph{Ran}(\overline{f}_\beta))\cap ACC_{0}(\emph{Ran}(\overline{f}_\beta))=\emptyset$
where the set $\mathfrak{\overline{T}}_\beta(\emph{Ran}(\overline{f}_\beta))$ is the $\overline{F}^\beta$-transformation of $\emph{ Ran}(\overline{f}_\beta)$,
$\overline{f}_\beta(x)=\frac{a_\beta+f_\beta(x)}{2}$ for any $x\in[s_\beta,t_\beta]$ and
\begin{equation*}
\overline{F}^\beta(x,y)=\left\{
  \begin{array}{ll}
    a_\beta+\frac{(b_\beta-a_\beta)}{2}F_\beta(\frac{2(x-a_\beta)}{b_\beta-a_\beta},\frac{2(y-a_\beta)}{b_\beta-a_\beta}) & \emph{if }(x,y)\in[a_{\beta},\frac{b_{\beta}+a_{\beta}}{2}]^{2}, \\
    \min\{x,y\} & \emph{if }(x,y)\in[a_{\beta},b_{\beta}]^{2}\setminus[a_{\beta},\frac{b_{\beta}+a_{\beta}}{2}]^{2};
  \end{array}
\right.
\end{equation*}

\item if $f(t_{\beta})> b_{\beta}$ then $\mathfrak{\underline{T}}_\beta(\emph{Ran}(f_\beta))\cap ACC_{0}(\emph{Ran}(f_\beta))=\emptyset$
where the set $\mathfrak{\underline{T}}_\beta(\emph{Ran}(f_\beta))$ is the $\underline{F}^\beta$-transformation of $\emph{ Ran}(f_\beta)$ and
\begin{equation*}
\underline{F}^\beta(x,y)=\left\{
  \begin{array}{ll}
    F^\beta(x,y) & \emph{if }(x,y)\in[a_{\beta},b_{\beta})^{2}, \\
    \min\{x,y\} & \emph{if }(x,y)\in[a_{\beta},b_{\beta}]^{2}\setminus[a_{\beta},b_{\beta})^{2};
  \end{array}
\right.
\end{equation*}
\end{enumerate}
\item if $t_{\beta_i}=s_{\beta_j}$ for some $\beta_i,\beta_j\in B_A^f$ and there exist two elements $u,v\in (s_{\beta_j},t_{\beta_j})$ such that $$F(f_{\beta_j}(u),f_{\beta_j}(v))\leq f(s_{\beta_j}^+),$$ then $f(x^-)\leq F(f(t_{\beta_i}),f(x))$ for all $x\leq t_{\beta_i}$;
\item $f(x^-)\leq F(f(1),f(x))$ for all $x\leq 1$.
\end{enumerate}
\end{theorem}

\section{Characterizations of classes of t-norms generated by \eqref{eq1}}

In 2002, Klement, Mesiar and Pap \cite{EP2002} had shown that all t-norms can be naturally divided into three classes, i.e., $\{T_M\}$, the class of ordinally irreducible t-norms and the class of t-norms which are different from $\{T_M\}$ and are not ordinally irreducible. In this section, we characterize each class of t-norms that are generated by \eqref{eq1}.

\begin{definition}[\cite{EP2000}]
\emph{A t-norm $T$ which only has a trivial ordinal sum representation (i.e., there is no ordinal sum representation of $T$ different from $(\langle0,1,T\rangle)$) is called ordinally irreducible.}
\end{definition}

\begin{remark}[ \cite{EP2002}]\label{rem6.1}
\emph{Based on Theorem \ref{theorem:4.1}, a t-norm can be seen as an ordinal sum of t-subnorms $(I_\alpha,*_\alpha)_{\alpha\in A}$ where $A\neq\emptyset$ is a totally ordered set and $I$ is a non-empty subinterval of $[0,1]$. Then all t-norms can be naturally divided into three classes as follows:}
\renewcommand{\labelenumi}{(\roman{enumi})}
\begin{enumerate}
\item \emph{$\{T_M\}$ where $T_M=\min\{x,y\}$ for all $x,y\in[0,1]$. In this case, $A$ is a singleton set, say $A=\{\alpha_1\}$, and $*_{\alpha_1}=\min$.}
\item \emph{The class of ordinally irreducible t-norms. In this case, $A$ is a singleton set, say $A=\{\alpha_1\}$, and $*_{\alpha_1}\neq\min$. For example,
\begin{equation*}
T^{nM}(x,y)=\left\{
  \begin{array}{ll}
    0& \hbox{if } x+y\leq 1, \\
     \min\{x,y\} &\hbox{otherwise}
  \end{array}
\right.
\end{equation*}
and
\begin{equation*}
T(x,y)=\begin{cases}
\max\{0.5,x+y-1\} & \hbox{if } x\in[0.5,1]^2,\\
\max\{0,x+y-1\} &  \mbox{otherwise}.\\
\end{cases}
\end{equation*}}
\item \emph{The class of t-norms which are different from $\{T_M\}$ and are not ordinally irreducible. In this case, $A$ has at least two elements or, equivalently, $([0,1],T)$ can be expressed as an ordinal sum of t-subnorm which contains at least one summand $(I,*)$ where the length of the interval $I$ is smaller than 1 and $*\neq\min$. For example, \begin{equation*}
T(x,y)=\left\{
  \begin{array}{ll}
    0.5+0.5\cdot \max\{2x+2y-3,0\} & \hbox{if }(x,y)\in[0.5,1]^{2}, \\
    \min\{x,y\} & \hbox{otherwise}
  \end{array}
\right.
\end{equation*}
and
\begin{equation*}
T(x,y)=\left\{
  \begin{array}{ll}
    0 & \hbox{if }(x,y)\in[0,0.5]^{2}, \\
    0.5+0.5\cdot (2x-1)(2y-1) & \hbox{if }(x,y)\in(0.5,1]^{2}, \\
    \min\{x,y\} & \hbox{otherwise.}
  \end{array}
\right.
\end{equation*}}
\end{enumerate}
\end{remark}

An interesting point is that a t-norm which is not ordinally irreducible can be obtained by \eqref{eq1} whenever one gives us a strictly increasing function $f:[0,1]\rightarrow[0,1]$ and an ordinally irreducible t-norm $F:[0,1]^2\rightarrow[0,1]$. For example, let $f(x)=0.5+0.5x$ for any $x\in[0,1]$ and $F=T^{nM}$. Then $T$ given by Eq.~\eqref{eq1} is $T_M$. Therefore, a natural question arises: what is a necessary and sufficient condition for the function $T$ given by Eq.~\eqref{eq1} being one of Remark \ref{rem6.1} (i),(ii) and (iii)? In the following, we give an answer of this problem.

\begin{theorem}\label{thm6.1}
Let $A\neq\emptyset$ be a totally ordered set, the t-norm $F=(\langle a_\alpha,b_\alpha, F_\alpha\rangle)_{\alpha\in A}$ be an ordinal sum of a family of t-subnorms $\{F_\alpha \mid \alpha\in A\}$, $f$ be a strictly increasing function and the set $\{f_\beta\mid f_\beta:[s_\beta,t_\beta]\rightarrow [a_\beta,b_\beta]\}_{\beta\in B_A^f}$ be a decomposition set of $f$. If a function $T$ is given by Eq.~\eqref{eq1}, then $T=T_M$ if and only if for any $y\in(0,1]$, $f(x^-)\leq F(f(x),f(y))$ for all $x\leq y$.
\end{theorem}
\begin{proof}We just prove the non-trivial part. Taking a $y\in(0,1]$, suppose that $f(x^-)\leq F(f(x),f(y))$ for all $x\leq y$. Then from $F(f(x),f(y))\leq \min\{f(x),f(y)\}$ we have $T(x,y)=x$ when $x\leq y$. Furthermore, by the monotonicity of $F$, we have $f(y^-)\leq F(f(y),f(y))\leq F(f(x),f(y))\leq f(y)$, and hence $T(x,y)=y$ when $y<x$. Therefore, $T(x,y)=\min\{x,y\}$ for any $x\in[0,1]$. Due to the arbitrariness of $y$, $T=T_M$.
\end{proof}

We need the following lemma.
\begin{lemma}[\cite{EP2002}]\label{lem6.1}
Let $T$ be a t-norm. Then the following are equivalent:
\renewcommand{\labelenumi}{(\roman{enumi})}
\begin{enumerate}
\item $T$ is ordinally irreducible.
\item For each $x\in(0,1)$ there exist $y,z\in[0,1]$ with $y<x<z$ and $T(y,z)<y$.
\end{enumerate}
\end{lemma}

\begin{theorem}\label{thm6.2}
Let $A\neq\emptyset$ be a totally ordered set, the t-norm $F=(\langle a_\alpha,b_\alpha, F_\alpha\rangle)_{\alpha\in A}$ be an ordinal sum of a family of t-subnorms $\{F_\alpha \mid \alpha\in A\}$, $f$ be a strictly increasing function and the set $\{f_\beta\mid f_\beta:[s_\beta,t_\beta]\rightarrow [a_\beta,b_\beta]\}_{\beta\in B_A^f}$ be a decomposition set of $f$. If a function $T$ is given by Eq.~\eqref{eq1}, then
$T$ is an ordinally irreducible t-norm if and only if $B_A^f$ is a singleton set, say $B_A^f=\{\beta\}$, satisfying $f(x)\in (a_\beta,b_\beta)$ for any $x\in(0,1)$ and all the following four statements are satisfied:
\renewcommand{\labelenumi}{(\roman{enumi})}
\begin{enumerate}
\item If $F^\beta$ is a t-norm on $[a_{\beta},b_{\beta}]^2$, then $\mathfrak{T}_\beta(\emph{Ran}(f_\beta))\cap ACC_{0}(\emph{Ran}(f_\beta))=\emptyset$
where $\mathfrak{T}_\beta(\emph{Ran}(f_\beta))$ is the $F^\beta$-transformation of $\emph{ Ran}(f_\beta)$;
\item If $F^\beta$ is a proper t-subnorm on $[a_{\beta},b_{\beta}]^2$, then
\begin{enumerate}
\item
if $f(t_{\beta})\leq b_{\beta}$ then $ \mathfrak{\overline{T}}_\beta(\emph{Ran}(\overline{f}_\beta))\cap ACC_{0}(\emph{Ran}(\overline{f}_\beta))=\emptyset$
where the set $\mathfrak{\overline{T}}_\beta(\emph{Ran}(\overline{f}_\beta))$ is the $\overline{F}^\beta$-transformation of $\emph{ Ran}(\overline{f}_\beta)$,
$\overline{f}_\beta(x)=\frac{a_\beta+f_\beta(x)}{2}$ for any $x\in[s_\beta,t_\beta]$ and
\begin{equation*}
\overline{F}^\beta(x,y)=\left\{
  \begin{array}{ll}
    a_\beta+\frac{(b_\beta-a_\beta)}{2}F_\beta(\frac{2(x-a_\beta)}{b_\beta-a_\beta},\frac{2(y-a_\beta)}{b_\beta-a_\beta}) & \emph{if }(x,y)\in[a_{\beta},\frac{b_{\beta}+a_{\beta}}{2}]^{2}, \\
    \min\{x,y\} & \emph{if }(x,y)\in[a_{\beta},b_{\beta}]^{2}\setminus[a_{\beta},\frac{b_{\beta}+a_{\beta}}{2}]^{2};
  \end{array}
\right.
\end{equation*}

\item if $f(t_{\beta})> b_{\beta}$ then $\mathfrak{\underline{T}}_\beta(\emph{Ran}(f_\beta))\cap ACC_{0}(\emph{Ran}(f_\beta))=\emptyset$
where the set $\mathfrak{\underline{T}}_\beta(\emph{Ran}(f_\beta))$ is the $\underline{F}^\beta$-transformation of $\emph{ Ran}(f_\beta)$ and
\begin{equation*}
\underline{F}^\beta(x,y)=\left\{
  \begin{array}{ll}
    F^\beta(x,y) & \emph{if }(x,y)\in[a_{\beta},b_{\beta})^{2}, \\
    \min\{x,y\} & \emph{if }(x,y)\in[a_{\beta},b_{\beta}]^{2}\setminus[a_{\beta},b_{\beta})^{2};
  \end{array}
\right.
\end{equation*}
\end{enumerate}
\item $f(x^-)\leq F(f(1),f(x))$ for all $x\leq 1$;
\item For each $x\in(0,1)$ there exist two elements $y,z\in[0,1]$ with $y<x<z$ and $F(f(y),f(z))<f(y^-)$.
\end{enumerate}
\end{theorem}
\begin{proof} Suppose that $T$ is an ordinally irreducible t-norm. Then $B_A^f$ is a singleton set, say $B_A^f=\{\beta\}$, satisfying $f(x)\in (a_\beta,b_\beta)$ for any $x\in(0,1)$. Otherwise, there exists an $x\in(0,1)$ such that $x<s_\beta$ or $x>t_\beta$. If $x<s_\beta$, then for all $y,z\in[0,1]$ with $y<s_\beta<z$, we get $f(y)<a_\beta$ and $f(z)>a_\beta$, implying $T(y,z)= f^{(-1)}(F(f(y),f(z)))=f^{(-1)}(\min\{F(f(y),f(z))\})=f^{(-1)}(f(y))=y$, violating $T$ is ordinally irreducible. Analogously, if $x>t_\beta$, then for all $y,z\in[0,1]$ with $y<t_\beta<z$ we get $T(y,z)=y$, thus violating $T$ is ordinally irreducible.

By Theorem \ref{thm5.2} and Lemma \ref{lem6.1}, the rest of proof is proved immediately.
\end{proof}

From Remark \ref{rem6.1} and Theorems \ref{thm6.1} and \ref{thm6.2}, we have the following theorem.
\begin{theorem}\label{thm6.3}
Let $A\neq\emptyset$ be a totally ordered set, the t-norm $F=(\langle a_\alpha,b_\alpha, F_\alpha\rangle)_{\alpha\in A}$ be an ordinal sum of a family of t-subnorms $\{F_\alpha \mid \alpha\in A\}$, $f$ be a strictly increasing function and the set $\{f_\beta\mid f_\beta:[s_\beta,t_\beta]\rightarrow [a_\beta,b_\beta]\}_{\beta\in B_A^f}$ be a decomposition set of $f$. If a function $T$ is given by Eq.~\eqref{eq1}, then
$T$ is a t-norm with a non-trivial ordinal sum representation and $T\neq T_M$ if and only if $\mathbf{either}$ that $B_A^f$ is a singleton set, say $B_A^f=\{\beta\}$, in this case, either $f(x)\in (a_\beta,b_\beta)$ for any $x\in(0,1)$ and there exists an $x\in(0,1)$ such that $f(y^-)\leq F(f(y),f(z))$ for all $y,z\in[0,1]$ with $y<x<z$ or there exists an $x\in(0,1)$ such that $f(x)\notin (a_\beta,b_\beta)$ $\mathbf{or}$ that $B_A^f$ has at least two elements, and all the following five statements are satisfied:
\renewcommand{\labelenumi}{(\roman{enumi})}
\begin{enumerate}
\item If $F^\beta$ is a t-norm on $[a_{\beta},b_{\beta}]^2$, then $\mathfrak{T}_\beta(\emph{Ran}(f_\beta))\cap ACC_{0}(\emph{Ran}(f_\beta))=\emptyset$
where $\mathfrak{T}_\beta(\emph{Ran}(f_\beta))$ is the $F^\beta$-transformation of $\emph{ Ran}(f_\beta)$;
\item If $F^\beta$ is a proper t-subnorm on $[a_{\beta},b_{\beta}]^2$, then
\begin{enumerate}
\item
if $f(t_{\beta})\leq b_{\beta}$ then $ \mathfrak{\overline{T}}_\beta(\emph{Ran}(\overline{f}_\beta))\cap ACC_{0}(\emph{Ran}(\overline{f}_\beta))=\emptyset$
where the set $\mathfrak{\overline{T}}_\beta(\emph{Ran}(\overline{f}_\beta))$ is the $\overline{F}^\beta$-transformation of $\emph{ Ran}(\overline{f}_\beta)$,
$\overline{f}_\beta(x)=\frac{a_\beta+f_\beta(x)}{2}$ for any $x\in[s_\beta,t_\beta]$ and
\begin{equation*}
\overline{F}^\beta(x,y)=\left\{
  \begin{array}{ll}
    a_\beta+\frac{(b_\beta-a_\beta)}{2}F_\beta(\frac{2(x-a_\beta)}{b_\beta-a_\beta},\frac{2(y-a_\beta)}{b_\beta-a_\beta}) & \emph{if }(x,y)\in[a_{\beta},\frac{b_{\beta}+a_{\beta}}{2}]^{2}, \\
    \min\{x,y\} & \emph{if }(x,y)\in[a_{\beta},b_{\beta}]^{2}\setminus[a_{\beta},\frac{b_{\beta}+a_{\beta}}{2}]^{2};
  \end{array}
\right.
\end{equation*}

\item if $f(t_{\beta})> b_{\beta}$ then $\mathfrak{\underline{T}}_\beta(\emph{Ran}(f_\beta))\cap ACC_{0}(\emph{Ran}(f_\beta))=\emptyset$
where the set $\mathfrak{\underline{T}}_\beta(\emph{Ran}(f_\beta))$ is the $\underline{F}^\beta$-transformation of $\emph{ Ran}(f_\beta)$ and
\begin{equation*}
\underline{F}^\beta(x,y)=\left\{
  \begin{array}{ll}
    F^\beta(x,y) & \emph{if }(x,y)\in[a_{\beta},b_{\beta})^{2}, \\
    \min\{x,y\} & \emph{if }(x,y)\in[a_{\beta},b_{\beta}]^{2}\setminus[a_{\beta},b_{\beta})^{2};
  \end{array}
\right.
\end{equation*}
\end{enumerate}
\item If $t_{\beta_i}=s_{\beta_j}$ for some $\beta_i,\beta_j\in B$ and there exist $u,v\in (s_{\beta_j},t_{\beta_j})$ such that $F(f_{\beta_j}(u),f_{\beta_j}(v))\leq f(s_{\beta_j}^+)$, then $f(x^-)\leq F(f(t_{\beta_i}),f(x))$ for all $x\leq t_{\beta_i}$;
\item There exists two elements $x,y\in(0,1)$ with $x\leq y$ such that $F(f(y),f(x))<f(x^-)$;
\item $f(x^-)\leq F(f(1),f(x))$ for all $x\leq 1$.
\end{enumerate}
\end{theorem}

\section{Conclusions}
In this article, we gave a necessary and sufficient condition for the function $T$ defined by Eq.(\ref{eq1}) being associative when $F$ is a t-norm (even a t-subnorm) and $f$ is a strictly increasing function (see Theorem \ref{thm3.2} and \ref{them5.1}). These results are applied for investigating the associativity of function $T$ given by Eq.~\eqref{eq1} when the t-norm $F$ is an ordinal sum of t-norms and an ordinal sum of t-subnorms in the sense of A. H. Clifford, respectively (see Theorems \ref{thm4.2} and \ref{thm5.2}). Comparing with \cite{EP2000,PV2005,PV2008,PV2013} which focus on t-norms generated by additive generators via Eq.~\eqref{eq:1.8}, the difference in essence is that their work holds over a semigroup $([0,\infty],+)$ with the usual addition operation ``+", and our results are true when the t-norm $F$ can be written as an ordinal sum of t-norms and t-subnorms in the sense of A. H. Clifford, respectively. Moreover, Proposition \ref{prop4.4} shows that the ordinal sum theorems for t-subnorms and t-norms are closely related to each other via a strictly increasing function. Theorem \ref{thm5.2} specifically tells us how to obtain a new t-norm through a t-norm (even a t-subnorm) and a strictly increasing function. Theorems \ref{thm6.1} and \ref{thm6.3} shows that a t-norm which is not ordinally irreducible can be obtained by \eqref{eq1} whenever one gives us a strictly increasing function $f:[0,1]\rightarrow[0,1]$ and an ordinally irreducible t-norm $F:[0,1]^2\rightarrow[0,1]$, and Theorem \ref{thm6.2} shows that a t-norm which is ordinally irreducible can be obtained by \eqref{eq1} whenever one gives us a strictly increasing function $f:[0,1]\rightarrow[0,1]$ and a t-norm $F:[0,1]^2\rightarrow[0,1]$ that is not ordinally irreducible. %It is worth pointing out that one can verify that all results obtained are hold when the function $f$ is strictly decreasing.
\section*{Statements and Declarations}
$\mathbf{Competing\quad Interests:}$ The authors declare that they have no known competing financial interests or personal relationships that could have appeared to
influence the work reported in this article.

\end{document}